\documentclass[notitlepage,leqno]{article}

\usepackage[latin1]{inputenc}
\usepackage{latexsym}
\usepackage{amssymb}
\usepackage[T1]{fontenc}	
\usepackage[english]{babel}
\usepackage{amsmath,amssymb}
\usepackage{amsfonts}
\usepackage{amscd}
\usepackage{empheq}
\usepackage{fancyvrb}
\usepackage{marvosym}
\usepackage{mathrsfs}
\usepackage{amsthm}
\usepackage{verbatim}
\usepackage{tikz}
\usetikzlibrary{chains}
\usepackage{color}
\usepackage{chngpage}
\usetikzlibrary{matrix,through,arrows}

\renewcommand{\a }{\alpha }

\newcommand{\D }{\Delta }

\newcommand{\e }{\varepsilon }

\newcommand{\rh }{\rho }

\newcommand{\Sig }{\Sigma}

\newcommand{\be}{\begin{equation}}
\newcommand{\ee}{\end{equation}}
\newenvironment{pf}{\noindent{\sc Proof}.\enspace}{\rule{2mm}{2mm}\medskip}
\newenvironment{pfn}{\noindent{\sc Proof}}{\rule{2mm}{2mm}\medskip}
\newtheorem{thm}{Theorem}[section]
\newtheorem{pro}[thm]{Proposition}
\newtheorem{definition}[thm]{Definition}
\newtheorem{lem}[thm]{Lemma}
\newtheorem{rem}[thm]{Remark}
\newtheorem{cor}[thm]{Corollary}

\newtheorem{conjecture}[thm]{Conjecture}
\newtheorem{example}[thm]{Example}

\renewcommand{\epsilon}{\varepsilon}
\renewcommand{\bullet}{\cdot}

\author{Alessandro Carlotto$^{a}$ and Andrea Malchiodi$^{b}$}
\title{WEIGHTED BARYCENTRIC SETS AND SINGULAR LIOUVILLE EQUATIONS ON COMPACT SURFACES}
\date{}
\begin{document}

\maketitle


\noindent $^a$ STANFORD UNIVERSITY, Department of Mathematics - Sloan Hall, 94305 Stanford, CA

\medskip

\noindent $^b$ SISSA - Via Bonomea 265, 34136 Trieste, ITALY

\

\begin{description}
\item[ Abstract.]
Given a closed two dimensional manifold, we prove a general existence result for a class of elliptic PDEs with exponential nonlinearities and negative Dirac deltas on the right-hand side, extending a theory recently obtained for the regular case. This is done by global methods: since the associated Euler functional is in general unbounded from below, we need to define a new model space, generalizing the so-called \textsl{space of formal barycenters} and characterizing (up to homotopy equivalence) its very low sublevels. As a result, the analytic problem is reduced to a topological one concerning the contractibility of this model space. To this aim, we prove a new functional inequality in the spirit of \cite{cli1} and then we employ a min-max scheme based on a cone-style construction, jointly with the blow-up analysis given in \cite{bmt} (after \cite{bt1} and \cite{bm}).
This study is motivated by abelian Chern-Simons theory in self-dual regime, or from the problem of prescribing the Gaussian curvature in presence of conical singularities (hence generalizing a problem raised by Kazdan and Warner in \cite{kw}).

\end{description}

\

\section{Introduction}

\medskip

\noindent
In the last five decades, much attention has been paid to partial differential equations arising in the context of Conformal Geometry.\newline
Some basic examples are obtained by the Laplace-Beltrami operator $\D_{g}$ on a compact Riemannian surface $(\Sigma,g)$: under a conformal change of metric, say $g\mapsto\widehat{g}=e^{2w}g$, it is well-known that the Gauss curvature transforms according to the law
\begin{displaymath}
K_{\widehat{g}}=e^{-2w}(-\D_{g}w+K_{g})
\end{displaymath}
and furthermore $\D_{\widehat{g}}=e^{-2w}\D_{g}$. Analytic methods allow, for instance, to prove the fundamental \textsl{Uniformization Theorem}, asserting that every compact surface carries a (conformal) metric of constant curvature.
One can ask a somehow dual question, namely whether a given $g$ such that $K_{g}$ is constant can be conformal to a metric with Gaussian curvature a given function $K_{\widehat{g}}$. This problem, named after Kazdan-Warner (see \cite{kw}) and also known as \textsl{Nirenberg problem} in the special case when $(\Sigma, g)$ is the standard sphere, is modeled by a Liouville type equation on our surface $(\Sigma,g)$
\begin{equation}
\label{reg}
-\Delta_{g}u=\rho\left(\frac{h(x)e^{2u}}{\int_{\Sigma}h(x)e^{2u}\,dV_{g}}-1\right)
\end{equation}
with $\rho$ a real parameter and $h:\Sigma\rightarrow\mathbb{R}$  a smooth function. However, one basic feature of this geometric problem is that such a $\rho=K_{g}$ is related to the topology of $\Sigma$ by means of the Gauss-Bonnet formula
\begin{displaymath}
\int_{\Sigma}K_{g}\,dV_{g}=2\pi\chi(\Sigma).
\end{displaymath}
Once we assume, without loss of generality, that $Vol_{g}(\Sigma)=1$, we have that this equation forces $K_{g}$ to attain values that are (some) integer multiples of $4\pi$: therefore, on Riemann surfaces, we say that $K_{g}$ is a \textsl{quantized parameter}.
\newline

We might generalize equation \eqref{reg} by adding to the right-hand side a finite linear combination of Dirac deltas and hence getting \textsl{singular Liouville equations}
\begin{equation}
\label{sing}
-\Delta_{g}u=\rho\left(\frac{h(x)e^{2u}}{\int_{\Sigma}h(x)e^{2u}\,dV_{g}}-1\right)-2\pi\sum_{i=1}^{m}\alpha_{j}(\delta_{p_{j}}-1)
\end{equation}
where $p_{j}\in\Sigma$ are some fixed points. This equation has a strong geometric flavor as well, since the extra terms can be
viewed as singularities in the Gauss curvature corresponding to a \textsl{local conical structure}, as can be justified via an extension of the Gauss-Bonnet formula (see \cite{tr}):
\begin{displaymath}
\int_{\Sigma}K_{g}^{reg}\,dV_{g}=2\pi\left[\chi\left(\Sigma\right)+\sum_{J}\alpha_{j}\right],
\end{displaymath}
with
\begin{equation}
\label{form}
K_{g}=\textrm{smooth function}- 2\pi\sum_{J \ \textrm{finite}}\alpha_{j}\delta_{p_{j}}, \quad \alpha_{j}\in\left(-1,0\right)
\end{equation}
the first summand in \eqref{form} being denoted above by $K_{g}^{reg}$.\newline

Equation \eqref{sing} also arises in the study of self-dual multivortices in the Electroweak Theory by Glashow-Salam-Weinberg \cite{lai}, where $u$ can be interpreted as the logarithm of the absolute value of the wave function and the points $p_{j}$'s are the \textsl{vortices}, where the wave function vanishes. This class of problems has proved to be relevant in other physical frameworks, such as the study of the statistical mechanics of point vortices in the mean field limit (\cite{kiess}, \cite{clmp1}, \cite{clmp2}) and the abelian Chern-Simons Theory, as discussed in \cite{tar4}.


The regular Liouville problem, under a positivity assumption for the function $h$, has a well-known variational structure: indeed \eqref{reg} is the Euler-Lagrange equation associated to the $C^{1}$ functional
\begin{equation}
\label{funct}
J_{\rho}(u)=\int_{\Sigma}\left|\nabla_{g}u\right|^{2}\,dV_{g}+2\rho\int_{\Sigma}u\,dV_{g}-\rho\log\int_{\Sigma}h(x)e^{2u}\,dV_{g}
\end{equation}
defined on the Sobolev space $H^{1}\left(\Sigma,g\right)$. The weak form of the \textsl{Moser-Trudinger inequality} (see \cite{mos})
\begin{equation}
\label{mtw}
\log\int_{\Sigma}e^{2(u-\overline{u})}\,dV_{g}\leq\frac{1}{4\pi}\int_{\Sigma}\left|\nabla_{g}u\right|^{2}\,dV_{g}+C_{\Sigma,g} \quad u\in H^{1}(\Sigma,g)
\end{equation}
guarantees that $J_{\rho}$ is well-defined on $H^{1}(\Sigma,g)$ for any value of $\rho\in\mathbb{R}.$
Moreover, $J_{\rho}$ is lower semi-continuous with respect to the weak topology of that space and so, since \eqref{mtw} gives \textsl{coercivity} of $J_{\rho}$ if $\rho<4\pi,$ we immediately get existence of \textsl{critical points} for this range of values and the corresponding solvability of \eqref{reg}. It is clear that such critical points are \textsl{global minima} for $J_{\rho}.$ Such a direct variational approach does not apply to the case $\rho\geq4\pi$ as can be seen by exhibiting explicit examples. Let $p\in\Sigma$ an arbitrary (but fixed) point and let $\lambda>0.$ We define a one-parameter family of \emph{bubbling} functions as follows:
\begin{equation}
\label{bubble}
\varphi_{\lambda,p}(y)=\log\left(\frac{\lambda}{1+\lambda^{2}d_{g}^{2}(p,y)}\right),
\end{equation}
where $d_{g}$ is the Riemannian distance defined on $\Sigma$ by means of $g$. These functions appear in different contexts, for instance in the study of the \textsl{Yamabe problem} (see \cite{lp} and references therein) and exhibit a \textsl{peaked} behavior as $\lambda$ goes to infinity, specifically $e^{2\varphi_{\lambda, p}}\rightharpoonup \pi\delta_{p}.$ Moreover, it is possible to analyze the asymptotics of the different terms in \eqref{funct} and get
\begin{displaymath}
\int_{\Sigma}\left|\nabla_{g}\varphi_{\lambda,p}\right|^{2}\,dV_{g}\simeq 8\pi\log \lambda;
\quad \quad
\int_{\Sigma}\varphi_{\lambda,p}\,dV_{g}\simeq -\log\lambda.
\end{displaymath}
This fact, taking into account that $\int_{\Sigma}h(\cdot)e^{2\varphi_{\lambda,p}\left(\cdot\right)}\,dV_{g}$ is bounded above and below by fixed positive constants (independent of $\lambda$), implies that $J_{\rho}(\varphi_{\lambda, p})\to-\infty$ as $\lambda\to+\infty$ when $\rho>4\pi$ and hence the claim. Therefore $J_{\rho}$ is \textsl{not coercive} for $\rho>4\pi$ and so there is no hope of finding \textsl{global minima} and we need to attack the problem by means of different techniques. In the related recent literature, two guidelines can be highlighted: on the one hand, topological methods relying on the degree theory by Leray-Schauder (see \cite{cldeg}), on the other purely variational methods based on an improvement of the Moser-Trudinger inequality \eqref{mtw}.
Considering this second line of research, a pretty exhaustive existence theorem has been presented in \cite{zind}. Let us give a short description of the conceptual path that has led to such a conclusion. \newline
Exploiting the variational structure described above, the basic idea is to study the topology of the sublevels of the functional $J_{\rho}$ in the non-coercive regime. If we are able to detect a change in such topology, we may hope then to infer existence results via deformation lemmas. In order to investigate the structure of very \textsl{low} sublevels of \eqref{funct}, we first need to consider how the constant on the right-hand side of \eqref{mtw} can be sharpened under extra assumptions on the involved function. Indeed, it was shown by Chen and Li in \cite{cli1} that the constant $1/(4\pi)$ can be improved whenever $u$ is in some sense concentrated in $l+1$ well-separated regions on $\Sigma$ (for positive $l$) getting for any $\e>0$
\begin{equation}
\label{mti}
\log\int_{\Sigma}e^{2(u-\overline{u})}\,dV_{g}\leq\frac{1}{4(l+1)\pi-\widetilde{\epsilon}}\int_{\Sigma}\left|\nabla_{g}u\right|^{2}\,dV_{g}+C
\end{equation}
where $C$ depends on $\e$ (see Lemma \ref{fund} for a precise statement). This result gives important information on the structure of sublevels of $J_{\rho}$ or, more precisely, on the concentration phenomena characterizing the functions belonging to sufficiently low sublevels. For instance, if $\rho\in\left(4\pi,8\pi\right)$ and $u$ belongs to a sufficiently low sublevel of $J_{\rho}$,  then this inequality implies that it has to be \textsl{conformally concentrated} on a single region, and this is precisely what happens for the bubbling functions. More generally, we come to the following concentration result:
\begin{pro}[\cite{cli1}, \cite{zind}]
\label{l2}
Assume $\rho\in\left(4k\pi, 4(k+1)\pi\right)$ for some $k\geq 1.$ Then, for any $\epsilon>0$ and $r>0$ there exists a sufficiently large positive constant $L:=L(\epsilon, r)$ such that for every $u\in H^{1}(\Sigma, g)$ with $J_{\rho}(u)\leq -L$ there are $k$ points on $\Sigma$ (say $p_{1,u},...,p_{k,u}$) so that
\begin{displaymath}
\frac{\int_{\Sigma\setminus\cup_{i=1}^{k}B_{r}(p_{i,u})}e^{2u}\,dV_{g}}{\int_{\Sigma}e^{2u}\,dV_{g}}<\e.
\end{displaymath}
\end{pro}
This gives a clear hint for the definition of a model space describing, up to homotopy equivalence, the global topology of the very low sublevels of $J_{\rh}$.
For any integer $k\geq 1$ we define the $k$-th set of formal barycenters of $\Sig$ as
\begin{displaymath}
\Sig_{k}:=\left\{\sum_{i=1}^{k}t_{i}\delta_{p_{i}}:\  \sum_{i=1}^{k}t_{i}=1,\ t_{i}\geq0,\ p_{i}\in\Sig \quad \forall i\in\left\{1,\ldots,k\right\}\right\}.
\end{displaymath}

It is naively clear that there is a natural identification $\Sigma_{1}\sim\Sigma$, and $\Sigma$ can be seen just as a special case of this construction. Each set $\Sigma_{k}$ is enriched with the weak topology as a subspace of the dual of $C^{1}(\Sigma, g)$. Such topology on $\Sigma_{k}$ is actually metrizable and the inherited structure is that of a \textsl{stratified set}, consisting of parts having different dimensions.
Moreover, we can exploit a well-known result asserting that if $\Sig$ is a compact surface with no boundary, then $\Sig_{k}$ is not contractible for any $k\geq1$ (see \cite{dm} for a sketch of the argument given in \cite{bc}): once we prove that $\Sig_{k}$ is homotopy equivalent to $J_{\rho}^{-L}=\left\{u\in H^{1}\left(\Sigma,g\right)| J_{\rho}(u)\leq -L\right\}$ (for $L\gg 1$), we get at once the non-contractibility of such low sublevels. When $\rho\in(4\pi,8\pi)$ the construction of similar homotopy maps is very easy: indeed the previous concentration result suggests that we can in fact \textsl{project} the functions belonging to the very low sublevels of $J_{\rho}$ to the manifold $\Sigma$ itself and, conversely, to any point of $\Sigma$ we can associate a corresponding bubbling function centered on that point and with a concentration parameter $\lambda$ determined in terms of depth of the sublevel (see \cite{djlw2}).
In the general case, we can map $\Sigma_{k}$ into $J_{\rho}^{-L}$ by defining for any $\sigma\in\Sigma_{k},$ $\sigma=\sum_{i=1}^{k}t_{i}\delta_{p_{i}},$ and $\lambda>0,$  the function $\varphi_{\lambda,\sigma}(y):\Sigma\to\mathbb{R}$ by
\begin{equation}
\label{one}
\varphi_{\lambda,\sigma}(y):=\log\sum_{i=1}^{k}t_{i}\left(\frac{\lambda}{1+\lambda^{2}d_{g}^{2}(p_{i},y)}\right)^{2}-\log\pi.
\end{equation}
These functions generalize the bubbles introduced above (see \eqref{bubble}). Moreover, it is possible to derive the desired approximation properties via a refined asymptotic analysis, as performed in \cite{mal1}, namely getting that for $\lambda\to+\infty$ one has that
$e^{\varphi_{\lambda,\sigma}}\rightharpoonup \sigma$ and $J_{\rho}(\varphi_{\lambda,\sigma})\rightarrow-\infty$ uniformly for $\sigma\in\Sigma_{k}.$

Conversely, we might define an application from low sublevels of $J_{\rho}$ to the approximation space $\Sigma_{k}$ and prove the homotopical triviality of the compositions with the operator $\Phi$ defined in terms of the functions in \eqref{one}. On the other hand, the topology of sufficiently \textsl{high} sublevels of $J_\rho$ turns out to be trivial. More precisely, we can state the following:

\begin{pro} [\cite{dm},\cite{mal2}]Suppose $\rho\in\left(4k\pi,4(k+1)\pi\right)$ for some $k\geq1.$ Then, there exist a threshold $L>0$ and a continuous projection $\Psi:J_{\rho}^{-L}\to\Sigma_{k}$ satisfying:
\begin{itemize}
\item{if $\left(u_{n}\right)_{n\in\mathbb{N}}\subseteq J_{\rho}^{-L}$ is such that $e^{2u_{n}}\rightharpoonup\sigma$ for some $\sigma\in\Sigma_{k},$ then $\Psi(u_{n})\rightharpoonup\sigma;$}
\item{for $\lambda$ sufficiently large the composition map $\Psi(\varphi_{\lambda,\cdot})$ is homotopic to the identity in $\Sigma_{k}$ and in addition $\Psi(\varphi_{\lambda,\cdot})\to Id|_{\Sigma_{k}}$ as $\lambda\to+\infty;$}
\item{for $\lambda$ sufficiently large the composition map $u\mapsto \varphi_{\lambda,\Psi(u)}$ is homotopic to the identity in $J_{\rho}^{-L}.$}
\end{itemize}
As a corollary, there exists $L>0$ such that $J_{\rho}^{-L}$ has the same homology as $\Sigma_{k}.$ Moreover, there exists $\overline{b}\in\mathbb{R}$ so large that $b\geq\overline{b}$ implies that the sublevel ${J_{\rho}}^{b}$ is a deformation retract of $H^{1}_{-}(\Sigma, g)$ (the subspace of $H^{1}(\Sigma,g)$ consisting of functions with null mean) and therefore has the homology of a point.
\end{pro}

When the \textsl{Palais-Smale condition} holds, it is well known that a difference of topology in the sublevels
of a functional yields existence of critical points, which is proved via the classical deformation lemma.
Unfortunately it is still an open problem whether the P-S condition is satisfied for $J_\rho$: however the problem
can be bypassed exploiting a method originally introduced by Struwe in \cite{s} and used for this functional also in \cite{djlw2}. M. Lucia in \cite{luc} obtained an alternative deformation lemma yielding existence of an approximating sequence $\left(w_{n}\right)$ of critical points of $J_{\rho_{n}}$ for some $\rho_{n}\to\rho$. This reduces all the problem to a blow-up analysis, which was in fact performed in \cite{bm} and later refined in \cite{ls1}, \cite{li}, \cite{cl1}
and \cite{cldeg}.
%
By means of all these tools, Djadli \cite{zind} was finally able to prove the solvability of \eqref{reg} for $\rho\in\left(4k\pi, 4\left(k+1\right)\pi\right)$.\newline

With respect to equation \eqref{sing}, much of the existing literature concerns asymptotic analysis or compactness of solutions (see for instance \cite{bt1}, \cite{bt2}, \cite{clsing}, \cite{tar1}, \cite{luzh2}), while relatively few results are available about existence. In this sense, some perturbative results are given in \cite{dem}, \cite{esp} and an approach via infinite-dimensional degree theory is under current investigation in \cite{cl3} (see also \cite{clsing}). Our goal here is to describe a large variational theory for this kind of equation, which mainly relies on improved Moser-Trudinger inequalities and min-max methods, well fitting with the study of the regular case.\newline

As a preliminary step, let us see how a variational structure can be recovered. To this aim, consider the Green's functions
of $\Delta_{g}$ with poles at $p_j$, namely the distributional solutions of
\begin{displaymath}
\Delta_{g}G_{p_{j}}=2\pi(\delta_{p_{j}}-1),
\end{displaymath}
which are well-known (see \cite{aub}) to exist and to be smooth away from the singularities.
%
%
Performing the substitution $\widetilde{u}:=u-\sum_{j=1}^{m}\alpha_{j}G_{p_{j}}$ \eqref{sing} transforms into
\begin{equation}\label{mod}
-\Delta_{g}\widetilde{u}+\rho=\rho\frac{\widetilde{h}(x)e^{2\widetilde{u}}}{\int_{\Sigma}\widetilde{h}(x)e^{2\widetilde{u}}\,dV_{g}} \quad \textrm{on}\ \Sigma,
\end{equation}
with $\widetilde{h}(x)=h(x)e^{2\sum_{j=1}^{m}\alpha_{j}G_{p_{j}}}.$ Due to the fact that $G_{p_{j}}\simeq \log d_{g}(x,p_{j})$ near $p_{j}$ we  find that
\begin{displaymath}\label{datm}
\widetilde{h}\geq 0; \quad \widetilde{h}(x)\simeq d_{g}(x,p_{j})^{2\alpha_{j}} \quad \textrm{near} \ p_{j}.
\end{displaymath}
As a result, \eqref{mod} is nothing but the Euler-Lagrange equation for the modified functional
\begin{equation}
\label{functm}
J_{\rho,\underline{\alpha}}(\widetilde{u})=\int_{\Sigma}\left|\nabla_{g}\widetilde{u}\right|^{2}\,dV_{g}+2\rho\int_{\Sigma}\widetilde{u}\,dV_{g}
-\rho\log\int_{\Sigma}\widetilde{h}(x)e^{2\widetilde{u}}\,dV_{g}, \qquad \widetilde{u}\in H^{1}(\Sigma,g)
\end{equation}
(where $\underline{\alpha}=\left(\alpha_{1},\ldots,\alpha_{m}\right)\in\mathbb{N}^{m}$) and so we can study existence questions by global variational methods.

Let us spend some words on the role played in equation \eqref{sing} by the parameters. In principle, we allow $\rho$ and also the $\alpha_{j}$'s to be real numbers. However, the change of variables we performed above motivates (due to obvious integrability conditions) the assumption $\alpha_{j}>-1$ for any $j\in\left\{1,\ldots,m\right\}$ and this will be always implicit in the sequel. However, this restriction is very natural with respect to the geometric problem since a cone at $p$ of angle $\theta \in (0, 2\pi)$ corresponds to a term of the form $-2\pi\alpha\delta_{p}$ in \eqref{sing}, with $\theta=2\pi(1+\alpha)$.

\

While the recent papers \cite{bdm} and in \cite{mr} (see also Corollary 6 in \cite{bt1}) treated existence for
positive $\alpha$'s, more interesting for the physical applications, here we consider
the case $\alpha_{j} \in (-1,0)$, which is geometrically more relevant. Some results in the coercive case were proved in (see \cite{tr})
via the following \textit{Troyanov's inequality},
valid for $\alpha>-1$, $p\in\Sigma$ and similar in spirit to \eqref{mtw}:
\begin{equation}
\label{mttr}
\log\int_{\Sigma}d_{g}(x,p)^{2\alpha}e^{2(u-\overline{u})}\,dV_{g}\leq\frac{1}{4\pi \min \left\{1,1+\alpha\right\}} \int_{\Sigma}\left|\nabla_{g}u\right|^{2}\,dV_{g}+C_{\alpha,\Sigma,g} \quad u\in H^{1}(\Sigma,g).
\end{equation}
Again, it is seen by defining suitable \textsl{singular bubbling functions} that the value of the
above constant is sharp.
Notice that when $\a < 0$ the constant is larger than $\frac{1}{4\pi}$, resulting in a worse loss
of coercivity of $J_{\rho,\underline{\alpha}}$ compared to the regular case: coercivity actually holds only when $\rho<4\pi\min_{j=1,\ldots,m}(1+\alpha_{j})$,  so the topology of low sublevels of the functionals needs to be studied with more refined strategies.\newline

In Section 2 of this paper, we prove a new general version of the Chen-Li inequality, which combines both \eqref{mtw} and \eqref{mttr} in a global setting, see Lemma \ref{clgen}. The inequality somehow \textit{localizes} the volume
control in terms of the Dirichlet energy: we get an \textit{amount} of $4 \pi$ near regular points, by \eqref{mtw},
and an amount of $4\pi(1+\a_j)$ near each singular point $p_j$, provided concentration of conformal volume occurs.
This result suggests the introduction of a \textit{weighted} model space for the singular problem,  $\Sigma_{\rho,\underline{\alpha}}$, which plays the same role as $\Sigma_{k}$ in the regular case.
\begin{definition}
Given a point $q\in\Sigma$ we define its weighted cardinality as follows:
\begin{displaymath}
\chi(q)=\left\{\begin{matrix} 1+\alpha_{j} & \mbox{if }q=p_{j} \ \mbox{for some} \ j=1,\ldots,m; \\ 1 & \mbox{otherwise}.
\end{matrix}\right.
\end{displaymath}
The cardinality of any finite set of (pairwise distinct) points on $\Sigma$ is obtained extending $\chi$ by additivity.
\end{definition}
This enables us to easily describe selection rules to determine admissibility conditions for specific barycentric configurations in dependence on the values of the $\alpha_{j}$'s and $\rho.$

\begin{definition}
Suppose all the parameters $\rho,\alpha_{1},\ldots,\alpha_{m}$ are fixed. We define the corresponding space of formal barycenters as follows
\begin{equation}
\Sigma_{\rho,\underline{\alpha}}=\left\{\sum_{q_{j}\in J}t_{j}\delta_{q_{j}}:\  \sum_{q_{j}\in J}t_{j}=1,\ t_{j}\geq0,\ q_{j}\in\Sigma \quad 4\pi\chi(J)<\rho\right\}.
\end{equation}
\end{definition}
Notice that since we are considering negative weights the topological structure of $\Sigma_{\rho,\underline{\alpha}}$
is in general richer than that of $\Sig_k$ and strongly depends on the values of the parameters $\rho$ and $\underline{\alpha}$. For instance, when $m=2, \alpha_{1}=\alpha_{2}=\alpha$ and $\rho>8\pi\left(1+\alpha\right),\ \rho>4\pi, \  \rho<4\pi(2+\alpha)$ we get that $\Sigma_{\rho, \underline{\alpha}}$ is roughly obtained gluing together a mirror image of $\Sigma$ and a linear handle joining the singular points $p_{1}$ and $p_{2}$.

This new phenomenon causes some difficulties in applying the procedure for the regular case described above,
relating low sublevels to barycentric sets. For example, it is much harder in our case to define continuous
projections from $J_{\rho,\underline{\alpha}}^{-L}$ ($L \gg 0$) onto $\Sigma_{\rho,\underline{\alpha}}$: this problem is addressed in Section 3. This requires a preliminary study of the topological properties of $\Sigma_{\rho,\underline{\alpha}}$ as a stratified set, mainly concerning how a partial ordering can be put on the class of substrata (Definition 3.1), the structure of the boundary of a given stratum (Lemmas 3.2 and 3.8) and the way different strata may intersect (Lemma 6.1). Moreover, the construction presented in \cite{dm} for auxiliary connecting homotopies that are needed to define the projector operators must be substantially modified in order to take care of the selection rules defined above: this is done in Lemma 3.5. The basic idea is that those constraints do not allow us to move Dirac masses in $\Sigma_{\rho,\underline{\alpha}}$ freely, since for instance moving a mass form a singular point to a regular one leads in general to a violation of the condition $4\pi\chi(J)<\rho$.

In Section 4 instead we embed an image of $\Sigma_{\rho,\underline{\alpha}}$ into low sublevels of $J_{\rho,\underline{\alpha}}^{-L}$ by constructing suitable test functions which, compared to those
in \eqref{one}, have to take into account the presence of singular points. This is done using a sort
of interpolation between \textit{regular bubbles} and \textit{singular bubbles} (which, we recall,
can be used to show the sharpness of \eqref{mtw} and \eqref{mttr} respectively) when their center
approaches some of the points $p_j$, see \eqref{gamma} and \eqref{multic}. This is a new feature
compared to \cite{bdm} and \cite{mr}, where the profiles of test functions were of uniform type.\newline

The constructions in Sections 3 and 4 allow us to derive some information on the topology of low
sublevels of $J_{\rho,\underline{\alpha}}$, and then to run min-max schemes as for the regular case.
The compactness results however have to be modified to take the singularities into account,
and rely on the results in \cite{bmt}. Precisely, they hold true
for $\rho \not\in \mathfrak{S}$, where $\mathfrak{S}$ is introduced in the definition below.

\begin{definition}
We say that $\overline{\rho}>0$ is a singular value for Problem \eqref{sing} if
\begin{equation}
\label{defsv}
\overline{\rho}=4\pi n+4\pi\sum_{i\in I}(1+\alpha_{i})
\end{equation}
for some $n\in\mathbb{N}$ and $I\subseteq\left\{1,\ldots,m\right\}$ (possibly empty) satisfying
$n+card\left(I\right)>0$. The set of singular values will be denoted by $\mathfrak{S}=\mathfrak{S}\left(\underline{\alpha}\right)$.
\end{definition}

We are now in position to state the main result of this paper, proved in Section 5,  which is the following.

\begin{thm}
\label{grande}
Suppose that the parameters $\underline{\alpha}\in\left(-1,0\right)^{m}$ and $\rho\in\mathbb{R}_{>0}\setminus\mathfrak{S}$
are such that the set $\Sigma_{\rho,\underline{\alpha}}$ is not contractible with respect to the topology of $C^{1}(\Sigma,g)^{\ast}$. Then Problem \eqref{sing} admits a solution $u$ such that $u=v+\sum_{j=1}^{m}\alpha_{j}G_{p_{j}}$ with $G_{p_{\cdot}}$ the Green functions defined above and $v\in C^{\gamma}(\Sigma,g)$, for any $\gamma\in\left[0,\gamma_{0}\right)$ with $\gamma_{0}\in\left(0,1\right)$, solving equation \eqref{mod}.
\end{thm}

In Section 6 we show by means of a large class of examples that the non-contractibility condition above is in fact \textsl{very frequently} satisfied, and we present a conjecture that aims at classifying the cases when $\Sigma_{\rho,\underline{\alpha}}$ is contractible in terms of simple algebraic relations involving $\rho$ and $\underline{\alpha}$. It has to be mentioned that after the review process of the present article was completed, we could actually obtain a proof of this conjecture, which will be the object of a forthcoming paper.  \newline

An announcement of the present results is given in the preliminary note \cite{cm}. \newline

\textsl{Notations}. Throughout this article, we will always deal with two sorts of distances: the Riemannian distance on the manifold $\left(\Sigma,g\right)$ is $d_{g}$, while the metric associated to the weak convergence in $\Sigma_{\rho,\underline{\alpha}}$ (defined in Section 3) is simply $d$ (refer to equation \eqref{distsup}).
The notation $B_{r}\left(p\right)$ stands for the metric ball on $\Sigma$ having center $p$ and radius $r$. 
We will always use the function space $H^{1}\left(\Sigma,g\right)$ and the symbol $\left\|\cdot\right\|$ stands for its seminorm
\begin{displaymath}
\left\|u\right\|=\left(\int_{\Sigma}\left|\nabla_{g}u\right|^{2}\,dV_{g}\right)^{1/2}.
\end{displaymath}
Since all the equations we are interested in are invariant by adding constants, we will normalize the functions conveniently so that either $\overline{u}=\frac{1}{Vol_{g}\left(\Sigma\right)}\int_{\Sigma}u\,dV_{g}$ vanishes, or $\int_{\Sigma}e^{2u}\,dV_{g}=1$ (regular case) and $\int_{\Sigma}\widetilde{h}e^{2u}\,dV_{g}=1$ (singular case). In the first case, by the Poincaré-Wirtinger inequality $\left\|\cdot\right\|$ is indeed a real norm and correspondingly $H^{1}_{-}\left(\Sigma,g\right)$ is the Hilbert space of null average functions belonging to $H^{1}\left(\Sigma,g\right)$.
Large positive constants are always denoted by $C$ and the exact value of $C$ is allowed to vary from formula to formula and also within the same line. When we want to stress the dependence on some parameter, we add subscripts to $C$, hence obtaining things like $C_{\delta}$, $C_{\varepsilon, r, \Sigma, g}$ and so on. Notice that also constants with subscripts are allowed to vary. Lastly, the cardinality of a set $I$ is denoted by $card(I)$, while $\chi(I)$ is the weighted cardinality defined in Section 2.\newline 

\textsl{Acknowledgments}. A. C. completed part of this work during his stays at SISSA in Trieste,
supported by the Scuola Normale Superiore and therefore wishes to express his gratitude to both
these institutions. A. M. has been  supported by the FIRB project {\em Analysis and Beyond} from
MiUR. Both authors are grateful to D. Ruiz for his suggestions on the constructions in Section 4.

\

\section{Improved inequalities}

\medskip

\noindent As anticipated in the introduction, the core of the variational approach to Problem \eqref{reg} is represented by an improvement of the Moser-Trudinger inequality first obtained by Chen and Li in \cite{cli1}: the constant $1/(4\pi)$ can be improved whenever $u$ is in some sense concentrated in well-separated regions on $\Sigma.$

\begin{lem} \label{fund}
Let $l$ be a positive integer, let $\Omega_{1},..,\Omega_{l+1}$ be disjoint subsets of $\Sigma$ satisfying a separation condition $d_{g}(\Omega_{i},\Omega_{j})>\delta_{0}$ for any $i\neq j$ and some $\delta_{0}>0$ and consider any $\gamma_{0}\in\left(0,\frac{1}{l+1}\right).$ Then, for any $\widetilde{\epsilon}>0,$ there exists a constant $C:=C(\Sigma, g, l, \delta_{0}, \gamma_{0}, \widetilde{\epsilon})$ such that
\begin{displaymath}
\log\int_{\Sigma}e^{2(u-\overline{u})}\,dV_{g}\leq\frac{1}{4(l+1)\pi-\widetilde{\epsilon}}\int_{\Sigma}\left|\nabla_{g}u\right|^{2}\,dV_{g}+C.
\end{displaymath}
for all functions $u\in H^{1}(\Sigma)$ satisfying
\begin{equation}
\label{cond}
\frac{\int_{\Omega_{i}}e^{2u}\,dV_{g}}{\int_{\Sigma}e^{2u}\,dV_{g}}\geq\gamma_{0}, \quad \forall i\in \left\{1,...,l+1\right\}.
\end{equation}
\end{lem}

The proof we are going to present here is significantly different from the one given by the authors in \cite{cli1} and is inspired on a \textsl{spectral decomposition} implemented by Djadli and Malchiodi in \cite{dm} for the Paneitz operator.
This is done because the same technique also fits the needs for the corresponding concentration inequalities in the singular case. Therefore we present it here both for the convenience of the reader and in order to make the proof of Lemma \ref{clgen}, regarding the singular case, more direct and conceptually clear.
\medskip

\begin{pf}
We only prove the result for $l=1,$ being the general case identical in the substance.\newline
It is possible to find two functions $k_{1}, k_{2}$ satisfying the following properties:
\begin{displaymath}
\left\{ \begin{array}{ll}
k_{i}(x)\in\left[0,1\right] & \textrm{for every $x\in\Sigma;$}\\
k_{i}(x)=1 & \textrm{for every $x\in\Omega_{i}, i=1,2;$}\\
k_{i}(x)=0 & \textrm{if $d(x,\Omega_{i})\geq\frac{\delta_{0}}{4};$} \\
\left\|k_{i}\right\|_{C^{2}(\Sigma, g)}\leq C_{\delta_{0},}
\end{array} \right.
\end{displaymath}
where $C_{\delta_{0}}$ is some positive constant just depending on $\delta_{0}$ ($C_{\delta_{0}}\sim 1/\delta_{0}^{2}$).

We first need some preparatory estimates, so fix a function $w\in H^{1}(\Sigma):$ without losing any generality, we can also assume that $\overline{w}=0$ and, by symmetry, that $\left\|k_{1}w\right\|\leq\left\|k_{2}w\right\|$. Using our hypothesis and \eqref{mtw}, we get

\begin{displaymath}
\int_{\Sigma}e^{2w}\,dV_{g}\leq\frac{1}{\gamma_{0}}\int_{\Omega_{1}}e^{2w}\,dV_{g}\leq\frac{1}{\gamma_{0}}\int_{\Sigma}e^{2k_{1}w}\,dV_{g}\leq\frac{C_{\Sigma, g}}{\gamma_{0}}\exp\left\{\frac{1}{4\pi}\left\|k_{1}w\right\|^{2}+\overline{k_{1}w}\right\}.
\end{displaymath}
Now, by construction $k_{1}w$ and $k_{2}w$ have well-separated supports and so in evaluating $\left\|(k_{1}+k_{2})w\right\|^{2}=\int_{\Sigma}\left|\nabla_{g}(k_{1}+k_{2})w\right|^{2}\,dV_{g}$ we do not have \textsl{mixed terms} and just get $\left\|(k_{1}+k_{2})w\right\|^{2}=\left\|k_{1}w\right\|^{2}+\left\|k_{2}w\right\|^{2}$ and consequently $\left\|k_{1}w\right\|^{2}\leq\frac{1}{2}\left\|(k_{1}+k_{2})w\right\|^{2}.$ Exploiting these two inequalities we get
\begin{equation}
\label{prov1}
\int_{\Sigma}e^{2w}\,dV_{g}\leq\frac{C_{\Sigma, g}}{\gamma_{0}}\exp\left\{\frac{1}{8\pi}\left\|(k_{1}+k_{2})w\right\|^{2}+\overline{k_{1}w}\right\}.
\end{equation}
Now, we need to work on these terms on the right-hand side of \eqref{prov1}. Concerning the average term, we use the classical \textsl{Young inequality} $ab\leq\epsilon a^{2}+\frac{1}{\epsilon}b^{2}$ (valid for \textsl{any} $\epsilon>0$) to get
\begin{displaymath}
\overline{k_{1}w}=\int_{\Sigma}k_{1}w\,dV_{g}\leq\int_{\Sigma}\left(\frac{k_{1}^{2}}{\epsilon}+\epsilon w^{2}\right)\,dV_{g}\leq\frac{1}{\epsilon}+\epsilon\left\|w\right\|_{2}^{2}.
\end{displaymath}
We then need to study the gradient terms, that can be handled separately. For instance
\begin{displaymath}
\int_{\Sigma}\left|\nabla_{g}(k_{1}w)\right|^{2}\,dV_{g}=\int_{\Sigma}\left|\left(\nabla_{g}k_{1}\right)w+k_{1}\left(\nabla_{g}w\right)\right|^{2}\,dV_{g}
\end{displaymath}
\begin{displaymath}
=\int_{\Sigma}\left|\nabla_{g}k_{1}\right|^{2}w^{2}\,dV_{g}+\int_{\Sigma}k_{1}^{2}\left|\nabla_{g}w\right|^{2}\,dV_{g}+2\int_{\Sigma}k_{1}w\nabla_{g}\left(k_{1}\right)\nabla_{g}\left(w\right)\,dV_{g}
\end{displaymath}
\begin{displaymath}
\leq C_{\delta_{0}}\int_{\textrm{supp$(k_{1})$}}w^{2}\,dV_{g}+\int_{\textrm{supp$(k_{1})$}}\left|\nabla_{g}w\right|^{2}\,dV_{g}+2\epsilon\int_{\textrm{supp$(k_{1})$}}\left|\nabla_{g}w\right|^{2}\,dV_{g}+2\frac{C_{\delta_{0}}}{\epsilon}\int_{\textrm{supp$(k_{1})$}}w^{2}\,dV_{g}
\end{displaymath}
again applying the Young inequality (for the same value of $\epsilon$).
Hence, this leads to
\begin{displaymath}
\left\|(k_{1}+k_{2})w\right\|^{2}\leq C_{\delta_{0}}\left(1+\frac{2}{\epsilon}\right)\left\|w\right\|_{2}^{2}+\left(1+2\epsilon\right)\left\|w\right\|^{2}
\end{displaymath}
and by just renaming $\epsilon\rightarrow 2\epsilon$ for the sake of clarity we come to the auxiliary estimate
\begin{equation}
\label{aux}
\int_{\Sigma}e^{2w}\,dV_{g}\leq \frac{C}{\gamma_{0}}\exp\left\{\frac{1}{8\pi}\left(1+\epsilon\right)\left\|w\right\|^{2}+C_{\delta_{0}, \epsilon}\left\|w\right\|_{2}^{2}\right\},
\end{equation}
(where $C:=C(\Sigma, g,\epsilon)$), that will be used in the sequel of this proof to conclude the argument. \newline

Now, assume a generic function $u$ is given and pick $\widetilde{C}_{\delta_{0},\epsilon}$ so that $C_{\delta_{0},\epsilon}/\widetilde{C}_{\delta_{0},\epsilon}<\epsilon.$ It is standard and well known (see, for instance, \cite{aub} as a reference) that the operator $-\Delta_{g}$ admits a complete system of eigenfunctions on $X=H_{-}^{1}(\Sigma, g)$ and call $\left(\lambda_{j}\right)_{j\in\mathbb{N}}$ its (monotone and increasing) sequence of eigenvalues. We can then decompose $u$ as follows:

\begin{displaymath}
\left\{ \begin{array}{ll}
u=u_{\textrm{low}}+u_{\textrm{high}};\\
u_{\textrm{low}}=\sum_{\lambda_{j}\leq \widetilde{C}_{\delta_{0},\epsilon}}\varphi_{j};\\
u_{\textrm{high}}=\sum_{\lambda_{j}>\widetilde{C}_{\delta_{0}, \epsilon}}\varphi_{j}; \\
-\Delta_{g}\varphi_{j}=\lambda_{j}\varphi_{j} \quad \forall j\in\mathbb{N}.
\end{array} \right.
\end{displaymath}
On the one hand a straightforward computation shows that
\begin{displaymath} \left\|u_{\textrm{high}}\right\|_{2}^{2}\leq\frac{\left\|u_{\textrm{high}}\right\|^{2}}{\widetilde{C}_{\delta_{0},\epsilon}},
\end{displaymath}
while on the other $u_{\textrm{low}}\in L^{\infty}(\Sigma, g)$ with $\left\|u_{\textrm{low}}\right\|_{\infty}\leq C_{\delta_{0},\epsilon}\left\|u_{\textrm{low}}\right\|_{2}.$ In fact, there is equivalence between these two norms because the inequality $\left\|\bullet\right\|_{2}\leq\left\|\bullet\right\|_{\infty}$ is trivial (recall that we are assuming $Vol_{g}(\Sigma)=1$), while the other comes from elliptic regularity referred to the generators $\varphi_{j}$ of the finite-dimensional vector space $V_{\delta_{0},\epsilon}:=\left\langle \varphi_{j}| \lambda_{j}\leq \widetilde{C}_{\delta_{0},\epsilon}\right\rangle.$ Consequently, we can exploit both these facts proceeding as follows
\begin{displaymath}
\int_{\Sigma}e^{2u}\,dV_{g}=\int_{\Sigma}e^{2(u_{\textrm{low}}+u_{\textrm{high}})}\,dV_{g}\leq e^{2\left\|u_{\textrm{low}}\right\|_{\infty}}\int_{\Sigma}e^{2u_{\textrm{high}}}\,dV_{g}
\end{displaymath}
\begin{displaymath}
\leq e^{2\left\|u_{\textrm{low}}\right\|_{\infty}}\frac{C}{\gamma_{0}e^{-2\left\|u_{\textrm{low}}\right\|_{\infty}}}\exp\left\{\frac{1}{8\pi}(1+\epsilon)\left\|u_{\textrm{high}}\right\|^{2}+C_{\delta_{0}, \epsilon}\left\|u_{\textrm{high}}\right\|_{2}^{2}\right\},
\end{displaymath}
since we can make use of \eqref{aux} because the function $u_{\textrm{high}}$ satisfies the condition \eqref{cond} with $\gamma_{0}':=\gamma_{0}e^{-2\left\|u_{\textrm{low}}\right\|_{\infty}}.$
Equivalently, we have come to
\begin{displaymath}
\log\int_{\Sigma}e^{2u}\,dV_{g}\leq C+4\left\|u_{\textrm{low}}\right\|_{\infty}+\left\{\frac{1}{8\pi}(1+\epsilon)\left\|u_{\textrm{high}}\right\|^{2}+C_{\delta_{0}, \epsilon}\left\|u_{\textrm{high}}\right\|_{2}^{2}\right\},
\end{displaymath}
but due to the \textsl{Poincar\'e-Wirtinger inequality} and the elementary inequality $\sqrt{a}\leq\epsilon a+1/\epsilon,$ this becomes
\begin{displaymath}
\log\int_{\Sigma}e^{2u}\,dV_{g}\leq C+4\epsilon\left\|u_{\textrm{low}}\right\|^{2}+\left\{\frac{1}{8\pi}(1+\epsilon)\left\|u_{\textrm{high}}\right\|^{2}+C_{\delta_{0}, \epsilon}\left\|u_{\textrm{high}}\right\|_{2}^{2}\right\}.
\end{displaymath}

Depending on our choice of $\widetilde{C}_{\delta_{0},\epsilon}$ the previous inequality is just
\begin{equation}
\label{fin}
\log\int_{\Sigma}e^{2u}\,dV_{g}\leq C+\left\{\frac{1}{8\pi}(1+4\epsilon)\left\|u\right\|^{2}\right\}
\end{equation} where again $C=C(\Sigma, g, \delta_{0}, \gamma_{0}, \varepsilon).$
By means of some elementary algebra on the right-hand side of \eqref{fin}, we can replace this result (obtained for any $\epsilon>0$) with the thesis \eqref{mti}.
\end{pf}

The first step of our study is then a similar improved inequality that is based on both \eqref{mtw} and \eqref{mttr} and is proved still by means of cut-off functions, but with some extra algebra.

\begin{lem}
\label{clgen}
Let $n\in\mathbb{N}$ and let $I\subseteq\left\{1,\ldots,m\right\}$ with $n+card\left(I\right)>0$, where $card(I)$ denotes the  cardinality of a set. Assume there exists $r>0$, $\delta_{0}>0$ and pairwise distinct points $\left\{q_{1},\ldots,q_{n}\right\}\subseteq\Sigma\setminus\left\{p_{1},\ldots,p_{m}\right\}$ such that:
\begin{itemize}
\item{for any couple
 $\left\{a,b\right\}\subseteq\left\{q_{1},\ldots,q_{n}\cup\left(\cup_{i\in I}p_{i}\right)\right\}$ with $a\neq b$ one has $dist_{g}(B_{r}\left(a\right),B_{r}\left(b\right))\geq4\delta_{0}$;}
\item{for any $a\in\left\{q_{1},\ldots,q_{m}\right\}$ one has $d_{g}(p_{i},B_{r}(a))\geq4\delta_{0}$ for any $i\in\left\{1,\ldots,m\right\}\setminus I$;}
\end{itemize}
and consider any $\gamma_{0}\in\left(0,\frac{1}{n+card\left(I\right)}\right)$.

Then, for any $\widetilde{\epsilon}>0$ there exists a constant $C:=C(\Sigma, g, n, I, r, \delta_{0}, \gamma_{0}, \widetilde{\epsilon})$ such that
\begin{equation}\label{eq:cl}
\log\int_{\Sigma}\widetilde{h}e^{2(u-\overline{u})}\,dV_{g}\leq\frac{1}{4\pi\left(n+\sum_{i\in I}(1+\alpha_{i})-\widetilde{\epsilon}\right)}{\int_{\Sigma}\left|\nabla_{g}u\right|^{2}\,dV_{g}}+C
\end{equation}
for all functions $u\in H^{1}(\Sigma)$ satisfying
\begin{displaymath}
\frac{\int_{B_{r}(a)}\widetilde{h}e^{2u}\,dV_{g}}{\int_{\Sigma}\widetilde{h}e^{2u}\,dV_{g}}\geq\gamma_{0}, \quad \forall \ a\in\left\{q_{1},\ldots,q_{n}\cup\left(\cup_{i\in I}p_{i}\right)\right\}.
\end{displaymath}
\end{lem}

\medskip

\begin{pfn}
To avoid repetitions, we limit ourselves to sketch the argument, since many details can be borrowed from the proof of Lemma \ref{fund}. Assume first for any ball we deal with we define a suitable cut-off function. Exploiting them as above, we come to the following partial estimates (that hold for any $\epsilon>0$ small enough):
\begin{itemize}
\item{If $a\in \left\{q_{1},\ldots,q_{n}\right\}$ then
\begin{equation}
\label{prima}
\int_{\Sigma}\widetilde{h}e^{2w}\,dV_{g}\leq C \exp\left[\frac{1}{4\pi}(1+2\epsilon)\left\|w\right\|^{2}_{B_{r+\delta_{0}}(a)}+C_{\delta_{0},\epsilon}\left(\left\|w\right\|_{2}^{2}\right)_{B_{r+\delta_{0}}(a)}\right];
\end{equation}
}
\item{If $a=p_{i}$ for some $i\in I$ then by \eqref{mttr}
\begin{equation}
\label{seconda}
\int_{\Sigma}\widetilde{h}e^{2w}\,dV_{g}\leq C \exp\left[\frac{1}{4\pi(1+\alpha_{i})}(1+2\epsilon)\left\|w\right\|^{2}_{B_{r+\delta_{0}}(a)}+C_{\delta_{0},\epsilon}\left(\left\|w\right\|_{2}^{2}\right)_{B_{r+\delta_{0}}(a)}\right].
\end{equation}
}
\end{itemize}
Assume now we raise each of the inequalities \eqref{prima} to the power $\lambda^{-1}>0$ and the $i$-th of the inequalities \eqref{seconda} to the power $\mu_{i}^{-1}>0$ with
\begin{equation}
\label{aritm}
\left\{ \begin{array}{ll}
\frac{n}{\lambda}+\sum_{i\in I}\frac{1}{\mu_{i}}=1\\
\frac{1}{\lambda}\sum_{j=1}^{n}\theta_{j}+\sum_{i\in I}\frac{\varphi_{i}}{\mu_{i}(1+\alpha_{i})}\leq \frac{\sum_{j=1}^{n}\theta_{j}+\sum_{i\in I}\varphi_{i}}{n+\sum_{i\in I}(1+\alpha_{i})}
\end{array} \right.
\end{equation}
with $\theta_{j}=\left\|w\right\|^{2}_{B_{r+\delta_{0}}(q_{j})}$ and $\varphi_{i}=\left\|w\right\|^{2}_{B_{r+\delta_{0}}(p_{i})}$. The algebraic problem \eqref{aritm} is indeed solvable by setting for instance
\begin{displaymath}
\lambda=n+\sum_{i\in I}(1+\alpha_{i}), \quad \quad \mu_{i}=\frac{\lambda}{1+\alpha_{i}},\ i\in I.
\end{displaymath}
Hence, by multiplication of all such inequalities we get the intermediate result (true for any $\epsilon>0$ sufficiently small):
\begin{equation}
\label{ausilia}
\log \int_{\Sigma}\widetilde{h}e^{2w}\,dV_{g}\leq C+\left[\frac{1}{4\pi\left(n+\sum_{i\in I}(1+\alpha_{i})\right)}(1+\epsilon)\left\|w\right\|^{2}+C_{\delta_{0},\epsilon}\left\|w\right\|_{2}^{2}\right].
\end{equation}
The strategy now is to follow almost verbatim the proof of Lemma \ref{fund} and so to exploit spectral analysis of $-\Delta_{g}$ on $H^{1}_{-}(\Sigma,g)$ to absorb the $L^{2}$ term into the Dirichlet energy. Once we have decomposed $u=u_{\textrm{low}}+u_{\textrm{high}},$ we just need to apply \eqref{ausilia} for $u_{\textrm{high}}$ to get the thesis.
\end{pfn}

\begin{rem}
It should be clear that the same arguments work also if we replace the balls centered at singular points with balls covering the singular points (i.e. centered at points near the singularities), provided we guarantee some separation condition as above. This remark is actually useful for the proof of Lemma \ref{l2sing} below.
\end{rem}

\

\section{Mapping sublevels of $J_{\rho,\underline{\alpha}}$ into $\Sigma_{\rho,\underline{\alpha}}$}

\medskip

Following the guide of the regular case, we were led to claim the structure of the very low sublevels of the functional $J_{\rho,\underline{\alpha}}$ according to the definition of $\Sigma_{\rho,\underline{\alpha}}$ given in Section 1. Thanks to the previous improved inequalities, we expect that $\Sigma_{\rho,\underline{\alpha}}$ is indeed homotopy equivalent to the very low sublevels of the functional $J_{\rho,\underline{\alpha}}$: we introduce here a non-trivial projection operator $\Psi: J_{\rho,\underline{\alpha}}^{-L}\to \Sigma_{\rho,\underline{\alpha}}$ (for some appropriate choice of $L$) and, in the next section, an embedding $\Phi:\Sigma_{\rho,\underline{\alpha}}\to J_{\rho,\underline{\alpha}}^{-L}$ so that the composition $\Psi\circ\Phi:\Sigma_{\rho,\underline{\alpha}}\hookleftarrow$ is (homotopy) equivalent to the identity on the same space. Although this fact does not imply the homotopy equivalence, it is however sufficient for our purposes.

The model for this construction is presented in article \cite{dm}, where something similar is done (in a \textsl{regular} setting) for the $Q$-curvature prescription problem. Our case is for some aspects much harder. This is due to two related problems: 1) the topology of $\Sigma_{\rho,\underline{\alpha}}$ is very complicated and depends drastically on the values of the parameters, 2) the definition of the projection is delicate, since it must respect the selection rules for the barycenters defined above. The role of these obstructions should be clear in the sequel.\newline

Again, it is worth mentioning that the construction we are going to present is quite easy if we consider some specific values of the parameters (see Section 6 for some examples), but becomes rather sophisticated if we want to work in full generality. \newline

Throughout this section, we will consider $\Sigma_{\rho,\underline{\alpha}}$ endowed with the weak topology corresponding to the duality with $C^{1}(\Sigma,g)$. It is easy to see that such topology is equivalently determined by the distance function

\begin{equation}
\label{distsup}
d:\Sigma_{\rho,\underline{\alpha}}\times \Sigma_{\rho,\underline{\alpha}} \to \mathbb{R}_{\geq 0}\ , \quad d(\sigma_{1},\sigma_{2})= \sup_{\left\|f\right\|_{C^{1}\left(\Sigma\right)}\leq 1}\left(\sigma_{1}-\sigma_{2},f\right).
\end{equation}
This will be a useful tool to perform some explicit computations. \newline

We need to start by introducing some notation.
For $k, l \in \mathbb{N}$ and a set of indices $\left\{i_{1},\ldots,i_{l}\right\}\subseteq \left\{1,\ldots,m\right\}$ satisfying the relation $4\pi \left[k+\sum_{1}^{l}\left(1+\alpha_{i_{j}}\right)\right]<\rho$ we define the set
\begin{displaymath}
\Sigma^{k,l}_{i_{1}\ldots i_{l}}=\left\{s_{1}\delta_{p_{i_{1}}}+\ldots+ s_{l}\delta_{p_{i_{l}}}+\sum_{j=1}^{k}t_{j}\delta_{q_{j}}\right\},
\end{displaymath}
where
\begin{itemize}
\item{$s_{j}\in\left[0,1\right]$ for any $j=1,\ldots,l$;}
\item{$t_{j}\in\left[0,1\right]$ for any $j=1,\ldots,k$;}
\item{$\sum_{j}s_{j}+\sum_{j}t_{j}=1$;}
\item{$q_{j}\in \Sigma,$ for any $j=1,\ldots,k$.}
\end{itemize}

\begin{definition}
\label{ordstr}
Given two triplets $\left(k_{1},l_{1},\iota_{1}\right)$ and $\left(k_{2},l_{2},\iota_{2}\right)$, we will write that $\Sigma^{k_{1},l_{1}}_{\iota_{1}}\preceq \Sigma^{k_{2},l_{2}}_{\iota_{2}}$ if $\Sigma^{k_{1},l_{1}}_{\iota_{1}}\subseteq \Sigma^{k_{2},l_{2}}_{\iota_{2}}$ or, equivalently, if $k_{2}\geq k_{1}$ and the set of indices represented by $\iota_{1}$ can be split into two subsets, say $\overline{\iota}_{1}$ and $\overline{\overline{\iota}}_{1}$, such that:
\begin{itemize}
\item{$\overline{\iota}_{1}\subseteq \iota_{2}$;}
\item{$card\left(\overline{\overline{\iota}}_{1}\right)\leq k_{2}-k_{1}$}.
\end{itemize}
\end{definition}

This definition will be commented and motivated below, after a more general introduction of the construction we are going to perform. \newline

For any choice of $\left(k,l,\iota\right)$ we simply write $d_{k,l,\iota}\left(\sigma\right)=d\left(\sigma, \Sigma^{k,l}_{\iota}\right),$ $\sigma\in \Sigma_{\rho,\underline{\alpha}}$.
Then, for $\epsilon>0$ we define
\begin{displaymath}
\Sigma^{k,l}_{\iota}\left(\epsilon\right)=\left\{\sigma\in \Sigma^{k,l}_{\iota} | \ d_{k',l',\iota'}\left(\sigma\right)>\epsilon \ \textrm{for any triplet} \ \left(k',l',\iota'\right) \ \textrm{such that} \ \Sigma^{k',l'}_{\iota'}\prec \Sigma^{k,l}_{\iota} \right\}.
\end{displaymath}
In case $\Sigma^{k,l}_{\iota}$ is such that no triplet $\left(k',l',\iota'\right)$ exists with $\Sigma^{k',l'}_{\iota'}\prec\Sigma^{k,l}_{\iota}$, then we just set

\begin{displaymath}
\Sigma^{k,l}_{\iota}\left(\epsilon\right):=\Sigma^{k,l}_{\iota}.
\end{displaymath}
Such triplets $\left(k,l,\iota\right)$ will be called \textsl{minimal} with respect to $\prec$.
\newline

Lastly, we need to introduce an important tool. For any $l$ points $x_{1},\ldots,x_{l}\in\Sigma$ which all lie in a small metric ball and $l$ non-negative numbers $\gamma_{1},\ldots,\gamma_{l}$, we consider convex combinations of the form $\sum_{i=1}^{l}\gamma_{i}x_{i},\  \sum_{i}\gamma_{i}=1.$ To do this, we make use of the embedding of $\Sigma$ into some Euclidean space $\mathbb{R}^{n}$ given by Whitney's theorem, take the corresponding convex combination of these points in $\mathbb{R}^{n}$ and project it into our embedded manifold identified with the manifold itself. If $d_{g}\left(x_{i},x_{j}\right)<\xi$ for any choice of $i,j$ with $\xi$ sufficiently small this operation is well defined and moreover $d_{g}\left(x_{i},\sum_{j}\gamma_{j}x_{j}\right)<2\xi$ for any $i=1,\ldots,l$.  Alternatively, in order to preserve distances, we could employ Nash's embedding theorem, but this is not strictly necessary.

We now give a first quantitative description of the set $\Sigma^{k,l}_{\iota}$.

\begin{lem}\label{regoleps} Let $\left(k,l,\iota\right)$ a non-minimal admissible triplet. Then for all $\epsilon>0$ sufficiently small the following property holds: if $\sigma\in \Sigma^{k,l}_{\iota}\left(\epsilon\right),\ \sigma=\sum_{i=1}^{k+l}c_{i}\delta_{z_{i}}$, then
\begin{displaymath}
c_{i}\geq\frac{\epsilon}{2};\quad d_{g}\left(z_{i},z_{j}\right)\geq\frac{\epsilon}{2}; \quad i,j=1,\ldots,k+l, \ i\neq j.
\end{displaymath}
\end{lem}

\begin{pfn}
We study the two inequalities separately. Assume by contradiction the first is false and so there exists an index $\underline{i}\in\left\{1,\ldots,k+l\right\}$ such that $c_{\underline{i}}<\frac{\epsilon}{2}$. Then for $\underline{\underline{i}}\in\left\{1,\ldots,k+l\right\}, \underline{\underline{i}}\neq\underline{i}$ we consider the element
\begin{displaymath}
\widehat{\sigma}=\left(c_{\underline{i}}+c_{\underline{\underline{i}}}\right)\delta_{z_{\underline{\underline{i}}}}+\sum_{i=1,\ldots,k+l,\ i\neq\underline{i},\underline{\underline{i}}} c_{i}\delta_{z_{i}}.
\end{displaymath}
Depending on $\underline{i}$, the element $\widehat{\sigma}$ will belong either to $\Sigma^{k-1,l}_{\star}$ or to $\Sigma^{k,l-1}_{\star}$ for some multi-index $\star$ but in any case to a stratum (say $\Sigma^{k',l'}_{\iota'}$) that \textsl{precedes} $\Sigma^{k,l}_{\iota}$ in the sense explained above (see Definition \ref{ordstr}). Moreover, for any function $f\in C^{1}(\Sigma)$ with $\left\|f\right\|_{C^{1}\left(\Sigma\right)}\leq 1$ one has clearly
\begin{displaymath}
\left|\left(\sigma-\widehat{\sigma},f\right)\right|\leq c_{\underline{i}}\left(\left|f(z_{\underline{i}})\right|+\left|f\left(z_{\underline{\underline{i}}}\right)\right|\right)\leq 2c_{\underline{i}}
\end{displaymath}
and hence, taking the supremum with respect to $f$, we deduce
\begin{displaymath}
\epsilon<d\left(\sigma,\Sigma^{k',l'}_{\iota'}\right)\leq d\left(\sigma,\widehat{\sigma}\right)\leq \sup_{f}\left|\left(\sigma-\widehat{\sigma},f\right)\right|\leq 2c_{\underline{i}}.
\end{displaymath}
This is a contradiction.\newline
Let us now turn to the second inequality. Assume that there are $z_{i},z_{j}\in\Sigma$ with $z_{i}\neq z_{j}$ and $d_{g}\left(z_{i},z_{j}\right)<\frac{\epsilon}{2}$. Observe that, without losing any generality, we can assume that either $z_{i}$ or $z_{j}$ is \textsl{not} a singular point, simply because we can reduce the problem to the case $\epsilon<\min_{p\neq p'}d_{g}\left(p,p'\right)$ where $p, p'$ are a couple of singular points, so $\left\{p,p'\right\}\subseteq\left\{p_{1},\ldots,p_{m}\right\}$. Therefore, we can define the element
\begin{displaymath}
\widehat{\sigma}=\left(c_{i}+c_{j}\right)\delta_{\frac{1}{2}z_{i}+\frac{1}{2}\delta_{z_{j}}}+\sum_{s=1,\ldots,k+l \ s\neq i,j}c_{s}\delta_{z_{s}}.
\end{displaymath}
Again, the element $\widehat{\sigma}$ belongs to a stratum $\Sigma^{k',l'}_{\iota'}$ that precedes $\Sigma^{k,l}_{\iota}$ and, for $\left\|f\right\|_{C^{1}\left(\Sigma\right)}\leq 1$ we obtain
\begin{displaymath}
\left|\left(\sigma-\widehat{\sigma},f\right)\right|\leq c_{i}\left|f(z_{i})-f\left(\frac{z_{i}+z_{j}}{2}\right)\right|+c_{j}\left|f\left(z_{j}\right)-f\left(\frac{z_{i}+z_{j}}{2}\right)\right|.
\end{displaymath}
Taking the supremum over $f$, we deduce
\begin{displaymath}
\epsilon<d\left(\sigma,\Sigma^{k',l'}_{\iota'}\right)=\sup_{f}\left|\left(\sigma-\widehat{\sigma},f\right)\right|\leq 2d(z_{i},z_{j})
\end{displaymath}
and this gives as well a contradiction, so the proof is complete.
\end{pfn}

\begin{cor}
\label{mnfld}
For any triplet $\left(k,l,\iota\right)$ such that the stratum $\Sigma^{k,l}_{\iota}$ is admissible and non-minimal and for any $\epsilon>0$ sufficiently small, the set $\Sigma^{k,l}_{\iota}\left(\epsilon\right)$ is a smooth open manifold of dimension $3k+l-1$.
\end{cor}

\begin{pfn}
The previous Lemma \ref{regoleps} guarantees that in case we consider $\Sigma^{k,l}_{\iota}\left(\epsilon\right)$ instead of $\Sigma^{k,l}_{\iota}$, then all the numbers $c_{i}$ are uniformly bounded away from zero and also the mutual distance between any two points $z_{i}, z_{j}$ is uniformly bounded from below. Therefore, recalling that the coefficients $c_{i}$ satisfy the constraint $\sum_{i=1}^{k+l}c_{1}=1$, each element of $\Sigma^{k,l}_{\iota}\left(\epsilon\right)$ can be smoothly parameterized by $2k$ coordinates locating the points $z_{i}$ and by $k+l-1$ coordinates identifying the numbers $c_{i}$.
\end{pfn}

\begin{rem}
The previous corollary involves only non-minimal strata, so one could at first wonder about minimal ones. But actually, one easily sees that they can only be of the form $\Sigma^{0,1}_{j}$ for some $j\in\left\{1,\ldots,m\right\}$. Each of these only consists of one point, so the topology of such strata is also clear.
\end{rem}

In the regular case the strata are \textsl{totally} ordered by their dimensions and in fact:
\begin{displaymath}
\Sigma^{1}\prec\Sigma^{2}\prec\ldots\prec\Sigma^{k-1}\prec\Sigma^{k}, \quad \rho\in\left(4k\pi,4\left(k+1\right)\pi\right).
\end{displaymath}

In the singular case the situation is less clear in general. Given $d\in\mathbb{N}$, we may have different strata having dimension $d$ and this is due to two possibilities:
\begin{enumerate}
\item{We may have couples $\Sigma^{k,l}_{\iota}, \Sigma^{k,l}_{\iota'}$ with $\iota\neq \iota'$;}
\item{We may have couples $\Sigma^{k_{1},l_{1}}_{\iota_{1}}, \Sigma^{k_{2},l_{2}}_{\iota_{2}}$ with $\left(k_{1},l_{1}\right)\neq\left(k_{2},l_{2}\right)$ but $3k_{1}+l_{1}=3k_{2}+l_{2}$.}
\end{enumerate}
It is easily seen, via explicit examples, that both phenomena may really occur.
\newline

We now want to move towards the construction of the projection operator. The central problem, recognized in \cite{dm}, is to obtain continuity when strata of different dimensions meet. To explain this, we may refer to a very elementary example. Assume we have a square (i.e. its boundary) in the plane. We may think of it as a stratified set with the four vertices of dimension 0 and the four edges of dimension 1.  Assume we want to define a projection from a $\delta$-neighborhood of this square to the square itself. This is easy if we consider the central portion of each side, but becomes non-trivial if we lie near a vertex. Indeed we can have a couple of points near a diagonal (and near such vertex) with arbitrarily small mutual distance and if we just patch together the projections along different sides, these points would be sent far. To avoid this, we need to proceed by increasing dimension of the strata and hence first project radially to the vertices and then (on the remaining portion of our $\delta$-neighborhood) orthogonally to the sides. However, if we want to obtain a continuous global map, these definitions have to match and so we need to determine four transition annuli in order to define homotopies between these two sorts of projections.
\newline

The hard point of the construction is to define suitable homotopies on \textsl{transition domains} and this is done by means of the following lemma, which is a variation on a result contained in \cite{dm}.

\begin{lem}
\label{omotcoll}
Let $\left(k,l,\iota\right)$ be a triplet such that $\Sigma^{k,l}_{\iota}$ is an admissible stratum and let $\epsilon>0$ be sufficiently small. Then there exists a number $\widehat{\epsilon}\ll\epsilon$, only depending on $\epsilon$ and $\left(k,l,\iota\right)$ and a map $U^{t}_{k,l,\iota}$ from the set
\begin{displaymath}
\Sigma^{\widehat{\epsilon},\epsilon}_{k,l,\iota}:=\left\{\sigma\in\Sigma_{\rho,\underline{\alpha}} | \ d(\sigma,\Sigma^{k,l}_{\iota}\left(\epsilon\right))<\widehat{\epsilon}\right\}
\end{displaymath}
into $\Sigma_{\rho,\underline{\alpha}}$ such that the following four properties hold true:
\begin{enumerate}
\item[(i)]{$U^{0}_{k,l,\iota}=Id$ and $U^{t}_{k,l,\iota}{|_{\Sigma^{k,l}_{\iota}}}=Id|_{\Sigma^{k,l}_{\iota}}$ for every $t\in\left[0,1\right]$;}
\item[(ii)]{$U^{1}_{k,l,\iota}\left(\sigma\right)\in \Sigma^{k,l}_{\iota}\left(\frac{\epsilon}{2}\right)$ for every $\sigma\in\Sigma^{\widehat{\epsilon},\epsilon}_{k,l,\iota}$;}
\item[(iii)]{$d\left(U^{0}_{k,l,\iota}\left(\sigma\right),U^{t}_{k,l,\iota}\left(\sigma\right)\right)\leq C_{k,l,\iota,\epsilon}\sqrt{\widehat{\epsilon}}$ for every $\sigma\in\Sigma^{\widehat{\epsilon},\epsilon}_{k,l,\iota}$ and $t\in\left[0,1\right]$;}
\item[(iv)]{If $\sigma\in\Sigma^{\widehat{\epsilon},\epsilon}_{k,l,\iota}\cap \Sigma^{k',l'}_{\iota'}$ for any stratum $\Sigma^{k',l'}_{\iota'}$ such that $\Sigma^{k,l}_{\iota}\preceq \Sigma^{k',l'}_{\iota'}$ then $U^{t}_{k,l,\iota}\left(\sigma\right)\in \Sigma^{k',l'}_{\iota'}$ for every $t\in\left[0,1\right]$.}
\end{enumerate}
\end{lem}

\begin{rem}
Some comments are in order. First of all, the idea of this lemma is that if an element $\sigma\in\Sigma_{\rho,\underline{\alpha}}$ is near the set $\Sigma^{k,l}_{\iota}\left(\epsilon\right)$, then it can be projected to $\Sigma^{k,l}_{\iota}\left(\frac{\epsilon}{2}\right)$. Secondly, it has to be remarked that the constant $C_{k,l,\iota,\epsilon}$ does not depend on $t$ and $\widehat{\epsilon}$.
Finally, notice that among the properties above, probably the most important is the last one, because it tells that the homotopy $U^{t}_{k,l,\iota}$ acts respecting the higher strata, which should be a pretty natural requirement. The idea of (partially) ordering the strata by dimension - which is probably the first one could think of - does not work because such a definition of $\preceq'$ would necessarily lead to a violation of property (iv) above. The reason for this violation is explained after the proof of Lemma \ref{omotcoll} by means of Remark \ref{contr}.
\end{rem}

\begin{pfn}
We have seen in Corollary \ref{mnfld} that $\Sigma^{k,l}_{\iota}\left(\frac{\epsilon}{4}\right)$ is a smooth (open) finite-dimensional manifold and so there exists a projection $P_{k,l,\iota}$ from the $\widehat{\epsilon}$-neighborhood in $\Sigma_{\rho,\underline{\alpha}}$ of $\Sigma^{k,l}_{\iota}\left(\epsilon\right)$ onto $\Sigma^{k,l}_{\iota}\left(\frac{\epsilon}{2}\right)$. Due to the non-trivial structure of $\Sigma^{k,l}_{\iota}$ (it is not convex) and to the fact that $C^{1}(\Sigma,g)^{\ast}$ is a Banach space, this is actually only a \textsl{quasi-projection}, in the sense that
\begin{equation}
\label{quasipr}
d\left(\sigma,P_{k,l,\iota}\left(\sigma\right)\right)\leq C_{k,l,\iota,\epsilon}d\left(\sigma,\Sigma^{k,l}_{\iota}\left(\epsilon\right)\right),  \quad  \sigma\in\Sigma^{\widehat{\epsilon},\epsilon}_{k,l,\iota}.
\end{equation}
This construction is done by means of the Implicit Function Theorem and a partition of unity.
To fix the notation, we just write
\begin{displaymath}
\sigma=\sum_{i}c_{i}\delta_{z_{i}}, \quad P_{k,l,\iota}\left(\sigma\right)=\sum_{i}d_{i}\delta_{w_{i}}.
\end{displaymath}
Notice that we choose not do distinguish (at the level of notation) between regular and singular points, but to use this uniform notation. Notice also that since we are assuming $P_{k,l,\iota}\left(\sigma\right)\in\Sigma^{k,l}_{\iota}\left(\frac{\epsilon}{2}\right)$, then by Lemma \ref{regoleps}
\begin{displaymath}
d_{i}\geq \frac{\epsilon}{4}, \quad d_{g}\left(w_{i},w_{j}\right)\geq\frac{\epsilon}{4} \ \forall  \ i\neq j.
\end{displaymath}
Recall also that both the coefficients $d_{i}$ and the points $w_{i}\in\Sigma$ depend continuously on $\sigma$.\newline
We are going to define the map $U^{t}_{k,l,\iota}$ in different steps and the idea is basically first to reduce the number of points we deal with (this is done by means of a map $\widetilde{T}^{t}_{k,l,\iota}$ and its normalization $T^{t}_{k,l,\iota}$) and then to move towards $P_{k,l,\iota}(\sigma)$ in two steps in order to avoid transitions on forbidden configurations (see the selection rules above), i.e. we do not want to go out of $\Sigma_{\rho,\underline{\alpha}}$. \newline
Hence we first define an auxiliary map $\widetilde{T}^{t}_{k,l,\iota}, \ \widetilde{T}^{t}_{k,l,\iota}\left(\sigma\right)=\sum \widetilde{d}_{i}\delta_{\widetilde{w}_{i}}$ which misses the normalization $\sum \widetilde{d}_{i}=1$ and then correct the error. This application basically neglects the points $z_{i}$ that are far from any of the $w_{j}$'s (by letting their coefficients gradually vanishing) and sends any of the other points (say $z_{\underline{i}}$) to a convex combination of the points $z_{i}$'s that lie in a suitably small neighborhood of the same $w_{\underline{j}}$. However, differently from the regular case, we have to be careful with the singular points. In fact, this strategy could possibly lead to replace a singular point with a regular point (the corresponding convex combination), which might not be allowed in $\Sigma_{\rho,\underline{\alpha}}$. This is the reason for the introduction of the blow-up function $\theta:\Sigma\to\left[0,+\infty\right]$ that is defined as follows:
\begin{displaymath}
\theta(x)=\prod_{j=1}^{m}\theta_{j}(x), \quad \theta_{j}(x)=\max\left\{1, \frac{\mu}{d_{g}\left(x,p_{j}\right)}\right\},
\end{displaymath}
for some scale parameter $\mu\ll\min_{p\neq p'}d_{g}\left(p,p'\right)$ where $p, p'$ are singular points on the manifold $\Sigma$.

In order to obtain continuity for $\widetilde{T}^{t}_{k,l,\iota}$, we need to introduce a small parameter $\eta\ll\epsilon$ (that will be fixed later and will be of order $\simeq C_{k,l,\iota,\epsilon}\sqrt{\widehat{\epsilon}}$) and define a smooth cut-off function $\omega_{\eta}$ satisfying the following properties
\begin{equation}
\label{cuteta}
\left\{\begin{array}{ll}
\omega_{\eta}(t)=1,  & \textrm{for} \ t\leq\frac{\eta}{16};\\
\omega_{\eta}(t)=0, &  \textrm{for} \ t\geq\frac{\eta}{8};\\
\omega_{\eta}(t)\in\left[0,1\right], & \textrm{for \ every}\ t\geq 0.
\end{array}\right.
\end{equation}

Hence, we set $\omega_{j,\eta}\left(x\right)=\omega_{\eta}\left(d_{g}\left(x,w_{j}\right)\right)$.
We also define the following quantities:
\begin{displaymath}
\mathcal{X}_{j}\left(\sigma\right)=\frac{1}{\sum_{z_{i}\in B_{\frac{\eta}{8}}\left(w_{j}\right)}\theta\left(z_{i}\right)\omega_{j,\eta}\left(z_{i}\right)c_{i}}\sum_{z_{i}\in B_{\frac{\eta}{8}}\left(w_{j}\right)}\theta\left(z_{i}\right)\omega_{j,\eta}\left(z_{i}\right)c_{i}z_{i},
\end{displaymath}
\begin{displaymath}
s_{i}\left(\sigma\right)=\frac{8}{\eta}d_{g}\left(z_{i},w_{j}\right)-1, \ \textrm{for}\ z_{i}\in B_{\frac{\eta}{4}}\left(w_{j}\right).
\end{displaymath}
Since for any couple of indices $\underline{j}\neq\underline{\underline{j}}$ we have $d_{g}\left(w_{\underline{j}},w_{\underline{\underline{j}}}\right)\geq\frac{\epsilon}{4}$ and since $\eta\ll\epsilon$, then for any $i$ there exists (at most) one point $w_{j}$ such that $z_{i}\in B_{\frac{\eta}{4}}\left(w_{j}\right)$. As a result, the number $s_{i}\left(\sigma\right)$ is well-defined. After all these preliminaries, we define the map $\widetilde{T}^{t}_{k,l,\iota}$ as
\begin{displaymath}
\widetilde{T}^{t}_{k,l,\iota}\left(\sigma\right)=\sum_{i=1}^{k+l}\widetilde{c}_{i}\left(\sigma,t\right)\delta_{\widetilde{z}_{i}\left(\sigma,t\right)},
\end{displaymath}
with
\begin{displaymath}
\widetilde{c}_{i}\left(\sigma,t\right)=\left\{\begin{array}{ll}
\left(1-t\right)c_{i},  & \textrm{if} \ z_{i}\in \Sigma\setminus B_{\frac{\eta}{8}\left(w_{j}\right)};\\
\left(\left(1-t\right)+t\omega_{j,\eta}\left(z_{i}\right)\right)c_{i}, & \textrm{if}\ z_{i}\in B_{\frac{\eta}{8}}\left(w_{j}\right)
\end{array}\right.
\end{displaymath}
and
\begin{displaymath}
\widetilde{z}_{i}\left(\sigma,t\right)=\left\{\begin{array}{ll}
z_{i},  & \textrm{if} \ z_{i}\in \Sigma\setminus B_{\frac{\eta}{4}\left(w_{j}\right)};\\
\left(1-t\right)z_{i}+t\left[s_{i}\left(\sigma\right)z_{i}+\left(1-s_{i}\left(\sigma\right)\right)\mathcal{X}_{j}\left(\sigma\right)\right],  & \textrm{if} \ z_{i}\in B_{\frac{\eta}{4}}\left(w_{j}\right)\setminus B_{\frac{\eta}{8}\left(w_{j}\right)};\\
\left(1-t\right)z_{i}+t\mathcal{X}_{j}\left(\sigma\right), & \textrm{if}\ z_{i}\in B_{\frac{\eta}{8}}\left(w_{j}\right).
\end{array}\right.
\end{displaymath}

\begin{figure}[h]
\centering
\includegraphics[width=0.7\textwidth]{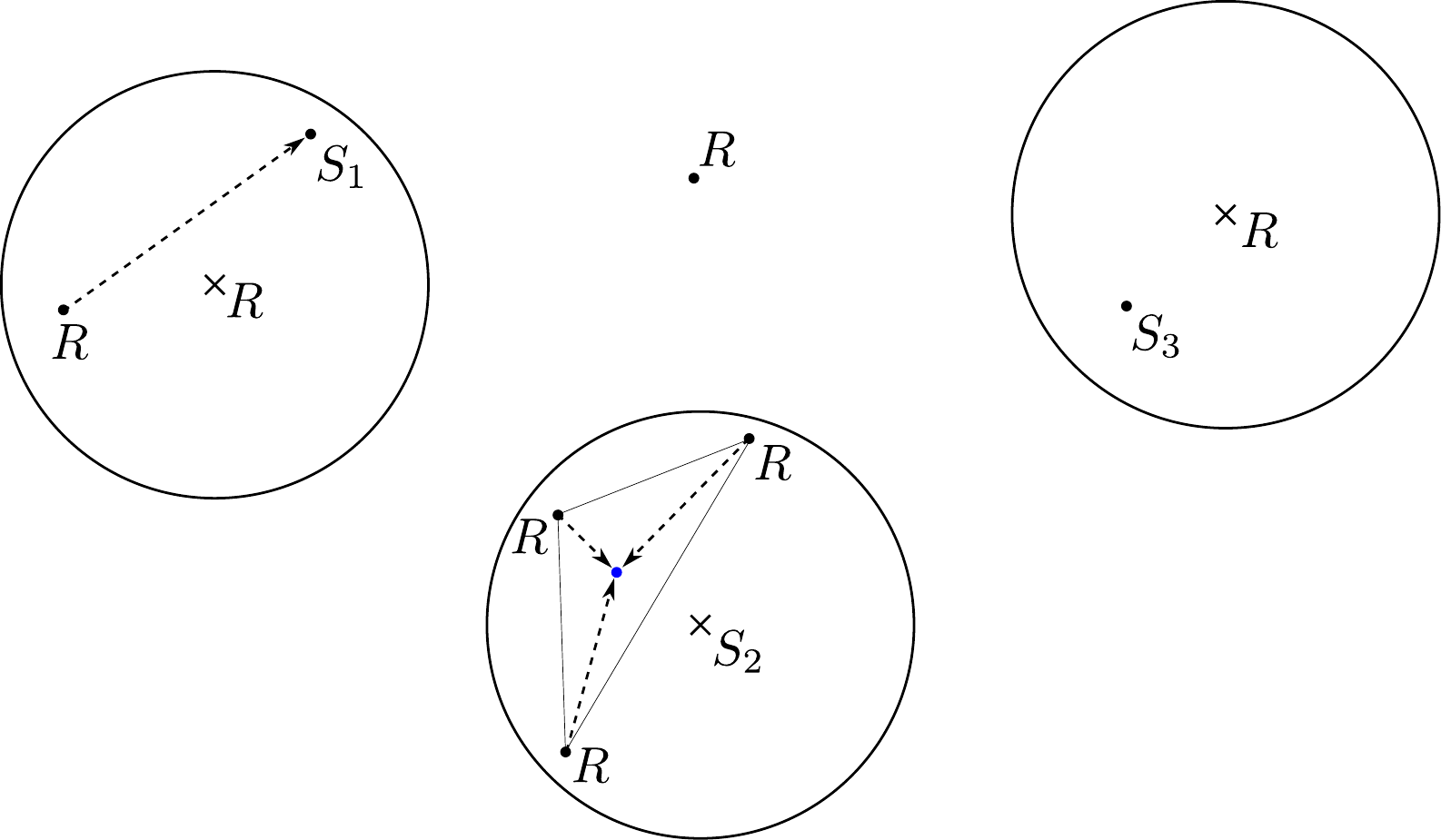}
\caption{The image through the projector $P_{2,1,2}$ of an element of the stratum $\Sigma^{5,2}_{13}$ and the line action of $\widetilde{T}^{t}_{2,1,2}$.}
\label{fig:balls}
\end{figure}

Now, the numbers $\widetilde{c}_{i}\left(\sigma,t\right)$ will in general miss the normalization condition $\sum_{i}\widetilde{c}_{i}=1$ and so we need to correct the map $\widetilde{T}^{t}_{k,l,\iota}$ defining
\begin{displaymath}
T^{t}_{k,l,\iota}\left(\sigma\right)=\frac{1}{\left(1-t\right)\widetilde{\mathcal{C}}\left(\sigma,0\right)+\sum\widetilde{\mathcal{C}}_{j}\left(\sigma,t\right)}\sum_{i=1}^{k+l}\widetilde{c}_{i}\left(\sigma,t\right)\delta_{\widetilde{z}_{i}\left(\sigma,t\right)},
\end{displaymath}
where
\begin{displaymath}
\widetilde{\mathcal{C}}_{j}\left(\sigma,t\right)=\sum_{z_{i}\in B_{\frac{\eta}{8}}\left(w_{j}\right)}\widetilde{c}_{i}\left(\sigma,t\right); \quad
\widetilde{\mathcal{C}}\left(\sigma,t\right)=1-\sum_{j}\widetilde{\mathcal{C}}_{j}\left(\sigma,t\right).
\end{displaymath}
One easily sees that the sum of all the coefficients is equal to 1 and that the map is well-defined and continuous in both $t$ and $\sigma$.
As a next step in our construction we need two more auxiliary maps. The first one is a homotopy  $H^{t}_{k,l,\iota}, t\in\left[0,1\right]$, that corrects the image of $T^{1}_{k,l,\iota}$ by sending the regular points among the $\widetilde{z}_{i}\left(\sigma,1\right)$'s to the corresponding image points $w_{j}$'s through $P_{k,l,\iota}$ and keeps the singular points still.
Lastly, we define a further correction homotopy $K_{k,l,\iota}$ so that each of the $\widetilde{z}_{i}\left(\sigma,1\right)$'s (and so the singular ones) is sent to its nearby image through $P_{k,l,\iota}$. The previous idea should be clear since the geometry of the set of points $\left(\cup_{i}\widetilde{z}_{i}\left(\sigma,1\right)\right)\cup\left(\cup_{j}w_{j}\right)$ is very simple and made of a finite number of couples (possibly singletons) contained in well-separated geodesic balls on $\Sigma$. Indeed, the definition of such homotopies $H_{k,l,\iota}$ and $K_{k,l,\iota}$ is elementary and we do not enter into details here.
We are now in position to complete our construction by setting
\begin{displaymath}
\widetilde{U}^{t}_{k,l,\iota}\left(\sigma\right)=\left\{\begin{array}{ll}
T^{3t}_{k,l,\iota},  & \textrm{for} \ t\in\left[0,\frac{1}{3}\right];\\
H^{3t-1}_{k,l,\iota}, & \textrm{for} \ t\in\left[\frac{1}{3},\frac{2}{3}\right];\\
K^{3t-2}_{k,l,\iota}, & \textrm{for}\ t\in\left[\frac{2}{3},1\right].
\end{array}\right.
\end{displaymath}
It is now needed to check the properties listed in the theorem. Among these, (i) is immediate, (iv) is easy and (ii) follows from (iii) (recall that we will finally make a smart choice of $\eta$ and $\widehat{\epsilon}$). So we just have to prove  property (iii) and it should be clear that we just need to verify it for the map $T^{t}_{k,l,\iota}$ since the action of both $H^{t}_{k,l,\iota}$ and $K^{t}_{k,l,\iota}$ is trivial and does not involve the coefficients.
\newline
This construction allows to adapt to our setting the estimates in \cite{dm}, that are reported here below for completeness. To begin, pick a smooth function $f$ such that
\begin{displaymath}
f\left(x\right)=\left\{\begin{array}{ll}
\frac{1}{2},  & \textrm{for} \ x\in \cup_{j}B_{\frac{\eta}{48}}\left(w_{j}\right);\\
\frac{1}{2}+\frac{\eta}{32}, & \textrm{for} \ x\in M\setminus \cup_{y}B_{\frac{\eta}{16}}\left(w_{j}\right);\\
\left\|f\right\|_{C^{1}\left(\Sigma\right)}\leq 1.
\end{array}\right.
\end{displaymath}
Since $\sigma\in \Sigma^{\widehat{\epsilon},\epsilon}_{k,l,\iota}$ and thanks to \eqref{quasipr} (that is $\left|\left(\sigma-P_{k,l,\iota}\left(\sigma\right),f\right)\right|\leq C_{k,l,\iota,\epsilon}\widehat{\epsilon}$) one has
\begin{equation}
\label{sommae}
\frac{\eta}{32} \sum_{z_{i}\in\Sigma\setminus \cup_{j}B_{\frac{\eta}{16}}\left(w_{j}\right)} c_{i}\leq \left(\sigma,f\right)-\left(P_{k,l,\iota},f\right)\leq C_{k,l,\iota,\epsilon}\widehat{\epsilon}
\end{equation}
because $\left(P_{k,l,\iota}\left(\sigma\right),f\right)=\sum_{j}d_{j}f\left(w_{j}\right)=1/2$ and
\begin{displaymath}
\left(\sigma,f\right)=\sum_{z_{i}\in \cup_{j}B_{\frac{\eta}{16}}\left(w_{j}\right)}c_{i}f\left(z_{i}\right)+\sum_{z_{i}\in M\setminus \cup_{j}B_{\frac{\eta}{16}}\left(w_{j}\right)}c_{i}f\left(z_{i}\right)
\end{displaymath}
\begin{displaymath}
\geq\frac{1}{2}\sum_{z_{i}\in\cup_{j}B_{\frac{\eta}{16}}\left(w_{j}\right)}c_{i}+\left(\frac{1}{2}+\frac{\eta}{32}\right)\sum_{z_{i}\in \Sigma\setminus \cup_{j}B_{\frac{\eta}{16}}\left(w_{j}\right)}c_{i}.
\end{displaymath}
The estimate \eqref{sommae} implies
\begin{displaymath}
\widetilde{\mathcal{C}}\left(\sigma,0\right)=\sum_{z_{i}\in\Sigma\setminus\cup_{j}B_{\frac{\eta}{8}}\left(w_{j}\right)}c_{i}\leq
\sum_{z_{i}\in\Sigma\setminus\cup_{j}B_{\frac{\eta}{16}}\left(w_{j}\right)}c_{i}\leq 32\frac{C_{k,l,\iota,\epsilon}\widehat{\epsilon}}{\eta}
\end{displaymath}
and also
\begin{displaymath}
\widetilde{\mathcal{C}}_{j}\left(\sigma,t\right)=\sum_{z_{i}\in B_{\frac{\eta}{8}}\left(w_{j}\right)\setminus B_{\frac{\eta}{16}}\left(w_{j}\right)}\left(\left(1-t\right)+t\omega_{j,\eta}\left(z_{i}\right)\right)c_{i}+\sum_{z_{i}\in B_{\frac{\eta}{16}}\left(w_{j}\right)}\left(\left(1-t\right)+t\omega_{j,\eta}\left(z_{i}\right)\right)c_{i}
\end{displaymath}
\begin{displaymath}
=\widetilde{\mathcal{A}}_{j}\left(\sigma,t\right)+\sum_{z_{i}\in B_{\frac{\eta}{16}}\left(w_{j}\right)}c_{i}, \ \textrm{where} \ \sum_{j}\widetilde{\mathcal{A}}_{j}\left(\sigma,t\right)\leq 32\frac{C_{k,l,\iota,\epsilon}\widehat{\epsilon}}{\eta}.
\end{displaymath}
Hence, exploiting the fact that by definition $\sum_{j}\widetilde{\mathcal{C}}_{j}\left(\sigma,0\right)+\widetilde{\mathcal{C}}\left(\sigma,0\right)=1$ or equivalently
\begin{displaymath}
\sum_{z_{i}\in B_{\frac{\eta}{16}}\left(w_{j}\right)}c_{i}+ \sum_{j}\widetilde{\mathcal{A}}_{j}\left(\sigma,0\right)+\sum_{z_{i}\in\Sigma\setminus \cup_{j}B_{\frac{\eta}{8}}\left(w_{j}\right)}c_{i}=1
\end{displaymath}
we deduce
\begin{displaymath}
\left|\sum_{j}\widetilde{\mathcal{C}}_{j}\left(\sigma,t\right)+\left(1-t\right)\widetilde{\mathcal{C}}\left(\sigma,0\right)-1\right|=
\end{displaymath}
\begin{displaymath}
\left|\left(\sum_{j}\left(\widetilde{\mathcal{A}}_{j}\left(\sigma,t\right)-\widetilde{\mathcal{A}}_{j}\left(\sigma,0\right)\right)\right)+\left(1-t\right)\sum_{z_{i}\in\Sigma\setminus \cup_{j}B_{\frac{\eta}{8}}\left(w_{j}\right)}c_{i}\right|
\end{displaymath}
\begin{displaymath}
\leq 64\frac{C_{k,l,\iota,\epsilon}\widehat{\epsilon}}{\eta}+32\frac{C_{k,l,\iota,\epsilon}\widehat{\epsilon}}{\eta}= 96\frac{C_{k,l,\iota,\epsilon}\widehat{\epsilon}}{\eta}.
\end{displaymath}
As a result, recalling the fact that $\widehat{\epsilon}$ will be chosen so small that $\frac{C_{k,l,\iota,\epsilon}\widehat{\epsilon}}{\eta}\ll 1$ we can use a Taylor expansion to conclude
\begin{displaymath}
\left|\frac{1}{\sum_{j}\widetilde{\mathcal{C}}_{j}\left(\sigma,t\right)+\left(1-t\right)\widetilde{\mathcal{C}}\left(\sigma,0\right)}-1\right|\leq 100\frac{C_{k,l,\iota,\epsilon}\widehat{\epsilon}}{\eta}.
\end{displaymath}
This is a very useful estimate because for an arbitrary function $f\in C^{1}\left(\Sigma\right)$ with $\left\|f\right\|_{C^{1}\left(\Sigma\right)}\leq 1$,
\begin{displaymath}
\left|\left(\sigma-T^{t}_{k,l,\iota}\left(\sigma\right),f\right)\right|\leq \left|\left(\sigma-\widetilde{T}^{t}_{k,l,\iota}\left(\sigma\right),f\right)\right|+\left|\left(\widetilde{T}^{t}_{k,l,\iota}\sigma-T^{t}_{k,l,\iota}\left(\sigma\right),f\right)\right|
\end{displaymath}
\begin{displaymath}
\leq \left|\left(\sigma-\widetilde{T}^{t}_{k,l,\iota}\left(\sigma\right),f\right)\right|+100 \frac{C_{k,l,\iota,\epsilon}\widehat{\epsilon}}{\eta}
\end{displaymath}
and so all we need to do is to evaluate the distance between $\sigma$ and $\widetilde{T}^{t}_{k,l,\iota}\left(\sigma\right)$. To this aim, observe that
\begin{displaymath}
\left|\left(\sigma-\widetilde{T}^{t}_{k,l\iota}\left(\sigma\right),f\right)\right| \leq \sum_{z_{i}\in\Sigma\setminus \cup_{j}B_{\frac{\eta}{4}}\left(w_{j}\right)}c_{i}
\end{displaymath}
\begin{displaymath}
\sum_{z_{i}\in \cup_{j}B_{\frac{\eta}{4}}\left(w_{j}\right)\setminus B_{\frac{\eta}{16}}\left(w_{j}\right)}\left|c_{i}f\left(z_{i}\right)-\widetilde{c}_{i}\left(\sigma,t\right)f\left(\widetilde{z}_{i}\left(\sigma,t\right)\right)\right|
\end{displaymath}
\begin{displaymath}
+\sum_{z_{i}\in B_{\frac{\eta}{16}}\left(w_{j}\right)}c_{i}d_{g}\left(z_{i},\widetilde{z}_{i}\left(\sigma,t\right)\right)
\end{displaymath}
(recall that we are working with test functions that are $1$-Lipschitz).
Now, the fact that $\eta$ is very small implies that
\begin{displaymath}
\left|c_{i}f\left(z_{i}\right)-\widetilde{c}_{i}\left(\sigma,t\right)f\left(\widetilde{z}_{i}\left(\sigma,t\right)\right)\right|\leq\left|c_{i}-\widetilde{c}_{i}\left(\sigma,t\right)\right|+\widetilde{c}_{i}\left(\sigma,t\right) d_{g}\left(z_{i},\widetilde{z}_{i}\left(\sigma,t\right)\right)\leq 2c_{i},
\end{displaymath}
and as a consequence
\begin{displaymath}
\left|\left(\sigma-\widetilde{T}^{t}_{k,l\iota}\left(\sigma\right),f\right)\right| \leq
2\sum_{\Sigma\setminus \cup_{j}B_{\frac{\eta}{16}}\left(w_{j}\right)}c_{i}+\sum_{j}\sum_{z_{i}\in B_{\frac{\eta}{16}}\left(w_{j}\right)}c_{i} d_{g}\left(z_{i},\mathcal{X}_{j}\left(\sigma\right)\right)
\end{displaymath}
\begin{displaymath}
\leq 64\frac{C_{k,l,\iota,\epsilon}\widehat{\epsilon}}{\eta}+\sum_{j}\sum_{z_{i}\in B_{\frac{\eta}{16}}\left(w_{j}\right)}c_{i}d_{g}\left(z_{i},\mathcal{X}_{j}\left(\sigma\right)\right).
\end{displaymath}
To estimate the last term we need a geometric argument based on our notion of convex combination on the abstract manifold $\Sigma$ (see above): we know that each point $z_{i}$ is shifted in the homotopy at most by $\eta/2$ and so exploiting the fact that $\sum_{i}c_{i}=1$ we conclude
\begin{displaymath}
\left|\left(\sigma-T^{t}_{k,l\iota}\left(\sigma\right),f\right)\right| \leq 164\frac{C_{k,l,\iota,\epsilon}\widehat{\epsilon}}{\eta}+\frac{\eta}{2}.
\end{displaymath}
This motivates the choice of $\eta=C_{k,,l,\iota,\epsilon}\sqrt{\widehat{\epsilon}}$ and that is the end of our proof.
\end{pfn}

It is now possible to give the anticipated motivation for our Definition \ref{ordstr}.

\begin{rem}
\label{contr}
Assume $m\geq 4$ and $\rho>4\pi,\ \rho>4\pi \sum_{i=1}^{4}\left(1+\alpha_{i}\right)$. This choice means that the space of formal barycenters $\Sigma_{\rho,\underline{\alpha}}$ contains, as special cases, the two strata $\Sigma^{1,0}$, having dimension 2 and $\Sigma^{0,4}_{1234}$ having dimension 3. Hence, by dimensional ordering $\Sigma^{1,0}\prec'\Sigma^{0,4}_{1234}$. Assume now we apply the previous Lemma \ref{omotcoll} to the stratum $\Sigma^{1,0}$: if property (iv) were true for $\preceq'$ then the homotopy $U^{t}_{1,0}$ should respect the higher-dimensional stratum $\Sigma^{0,4}_{1234}$ in the sense that for any $t\in\left[0,1\right]$ it should take values in $\Sigma^{0,4}_{1234}$ whenever applied to a point of the stratum itself.
Unfortunately, this is in contradiction with property (ii) because we require $U^{1}_{k,l,\iota}\in\Sigma^{k,l}_{\iota}\left(\frac{\epsilon}{2}\right)$ and clearly $\Sigma^{1,0}\left(\frac{\epsilon}{2}\right)\subseteq\Sigma^{1,0}\setminus\left(\Sigma^{0,1}_{1}\cup\Sigma^{0,1}_{2}\cup\Sigma^{0,1}_{3}\cup\Sigma^{0,1}_{4}\right)$.
\end{rem}

The basic idea to go further is the following: if for some element of $C^{1}\left(\Sigma,g\right)^{*}$ both the projections $P_{k,l,\iota}$ and $P_{k',l'.\iota'}$ are defined, with $\Sigma^{k,l}_{\iota}\prec\Sigma^{k',l'}_{\iota'}$, then we can consider the composition $U^{t}_{k,l\iota}\circ P_{k',l',\iota'}$ to get an homotopy between $P_{k,l,\iota}$ and $P_{k',l'.\iota'}$. In other terms $U^{t}_{k,l,\iota}$ is the \textsl{transition operator} we were looking for.\newline

We need two more technical lemmas. 

\begin{lem}
\label{esistpr}
For any $\epsilon$ sufficiently small, there exists $\widehat{\epsilon}$ such that it is possible to define a continuous projection from the set
\begin{displaymath}
\left\{f\in L^{1}\left(\Sigma\right)|\  f\geq0, \ \int_{\Sigma}f\,dV_{g}=1, \ d\left(f,\Sigma^{k,l}_{\iota}\left(\epsilon\right)\right)<\widehat{\epsilon}\right\}
\end{displaymath}
into $\Sigma^{k,l}_{\iota}\left(\frac{\epsilon}{2}\right)$.
\end{lem}

This first one is based on the fact that all the strata $\Sigma^{k,l}_{\iota}$ are finite-dimensional (see Corollary \ref{mnfld}).
The second concerns the intersections of different strata and tells that transition homotopies are needed only for couples of strata $\Sigma^{k,l}_{\iota}$ and $\Sigma^{k',l'}_{\iota'}$ such that $\Sigma^{k,l}_{\iota}\prec\Sigma^{k',l'}_{\iota'}$ and not, for instance, whenever $\textrm{dim}\left(\Sigma^{k,l}_{\iota}\right)<\textrm{dim}\left(\Sigma^{k',l'}_{\iota'}\right)$.

\begin{lem}
\label{inters}
Let $\Sigma^{k_{1},l_{1}}_{\iota_{1}}$ and $\Sigma^{k_{2},l_{2}}_{\iota_{2}}$ be strata that are included in $\Sigma_{\rho,\underline{\alpha}}$ for some fixed admissible values of $\underline{\alpha}$ and $\rho$. Then $\Sigma^{k_{1},l_{1}}_{\iota_{1}}\cap \Sigma^{k_{2},l_{2}}_{\iota_{2}}$ equals the union of all and only the strata that are contained both in $\Sigma^{k_{1},l_{1}}_{\iota_{1}}$ and $\Sigma^{k_{2},l_{2}}_{\iota_{2}}$, that are those $\Sigma^{k,l}_{\iota}$ such that $\Sigma^{k,l}_{\iota}\preceq\Sigma^{k_{1},l_{1}}_{\iota_{1}}$ and $\Sigma^{k,l}_{\iota}\preceq\Sigma^{k_{2},l_{2}}_{\iota_{2}}$.
\end{lem}

The proof of this result is straightforward, but still we decided to present this fact as a separate lemma in order to emphasize how easily intersections and boundary relations among strata can be treated by simply referring to the triplets $\left(k,l,\iota\right)$. 

There is still one missing tool which is needed for the construction of the global projection $\Psi$ from a suitable sublevel of $J_{\rho,\underline{\alpha}}$ to $\Sigma_{\rho,\underline{\alpha}}$. Indeed, in the case of the regular problem the inequality given by Lemma \ref{fund} is used to prove Proposition \ref{l2}, concerning the concentration phenomena characterizing functions belonging to very low sublevels of $J_{\rho}$, and this is clearly a preliminary step for defining a projector onto $\Sigma_{k}$. To that aim, the following lemma is needed, which works as a criterion implying the condition requested for applying Lemma \ref{fund}.

\begin{lem}[\cite{dm}]
\label{suff}
Let $l$ be a positive integer and consider a couple of positive numbers $\epsilon$ and $r$. Then, for any non-negative $f\in L^{1}(\Sigma)$ (normalized to $\left\|f\right\|_{1}=1$) satisfying
\begin{displaymath}
\int_{\cup_{j=1}^{l}B_{r}(p_{j})}f\,dV_{g}<1-\epsilon \quad \textrm{for every $l-$tuple} \ p_{1},...,p_{l}\in\Sigma
\end{displaymath}
there exist $\overline{\epsilon}<\epsilon,$ $\overline{r}>0$ and points $\overline{p}_{1},...,\overline{p}_{l+1}$ all depending on $(\Sigma, g)$ and $\epsilon, r, l$ and, just in the case of these points, also on $f$ such that
\begin{displaymath}
\int_{B_{\overline{r}}(\overline{p}_{1})}f\,dV_{g}\geq\overline{\epsilon},...,\int_{B_{\overline{r}}(\overline{p}_{l+1})}f\,dV_{g}\geq\overline{\epsilon}; \quad  B_{2\overline{r}}(\overline{p}_{i})\cap B_{2\overline{r}}(\overline{p}_{j})=\emptyset \quad \textrm{for} \ {i\neq j}.
\end{displaymath}
\end{lem}

In the singular case, the very same strategy does not apply, but we can nevertheless use this lemma and the improved inequality \eqref{eq:cl}, to get, by a tedious argument that we omit, the following result.

\begin{lem}
\label{l2sing}
For arbitrarily small $\epsilon>0$ and $r>0$ there exists a sufficiently large constant $L:=L(\epsilon,r)$ such that for every $u\in H^{1}\left(\Sigma,g\right)$ with $J_{\rho,\underline{\alpha}}\left(u\right)\leq -L$ there is a stratum $\Sigma^{k,l}_{i_{1}\ldots i_{l}}\subset\Sigma_{\rho,\underline{\alpha}}$ such that
\begin{displaymath}
\frac{\int_{B_{r}(a_{i_{1}})\cup\ldots\cup B_{r}(a_{i_{l}})\cup B_{r}(b_{1})\cup\ldots\cup B_{r}(b_{k})}\widetilde{h}e^{2u}\,dV_{g}}{\int_{\Sigma}\widetilde{h}e^{2u}\,dV_{g}}\geq 1-\epsilon,
\end{displaymath}
for some points $a_{i_{s}}\in B_{2r}(p_{i_{s}})$ ($s=1,\ldots,l$) and $b_{1},\ldots,b_{k}$ satisfying
\begin{equation}
\label{separa}
\min_{j=1,\ldots,k}\min_{i=1,\ldots,m} d_{g}(b_{j},p_{i})\geq2r.
\end{equation}
\end{lem}

As a consequence, we may come to the conclusion of this section.

\begin{lem}
\label{proi}
For any choice of $\rho$ and $\underline{\alpha}$ according to the restriction of Problem \eqref{sing}, there exists a large $\widehat{L}>0$ and a continuous map from $J_{\rho,\underline{\alpha}}^{-\widehat{L}}$ into $\Sigma_{\rho,\underline{\alpha}}$.
\end{lem}

\begin{figure}[h]
\centering
\includegraphics[width=0.8\textwidth]{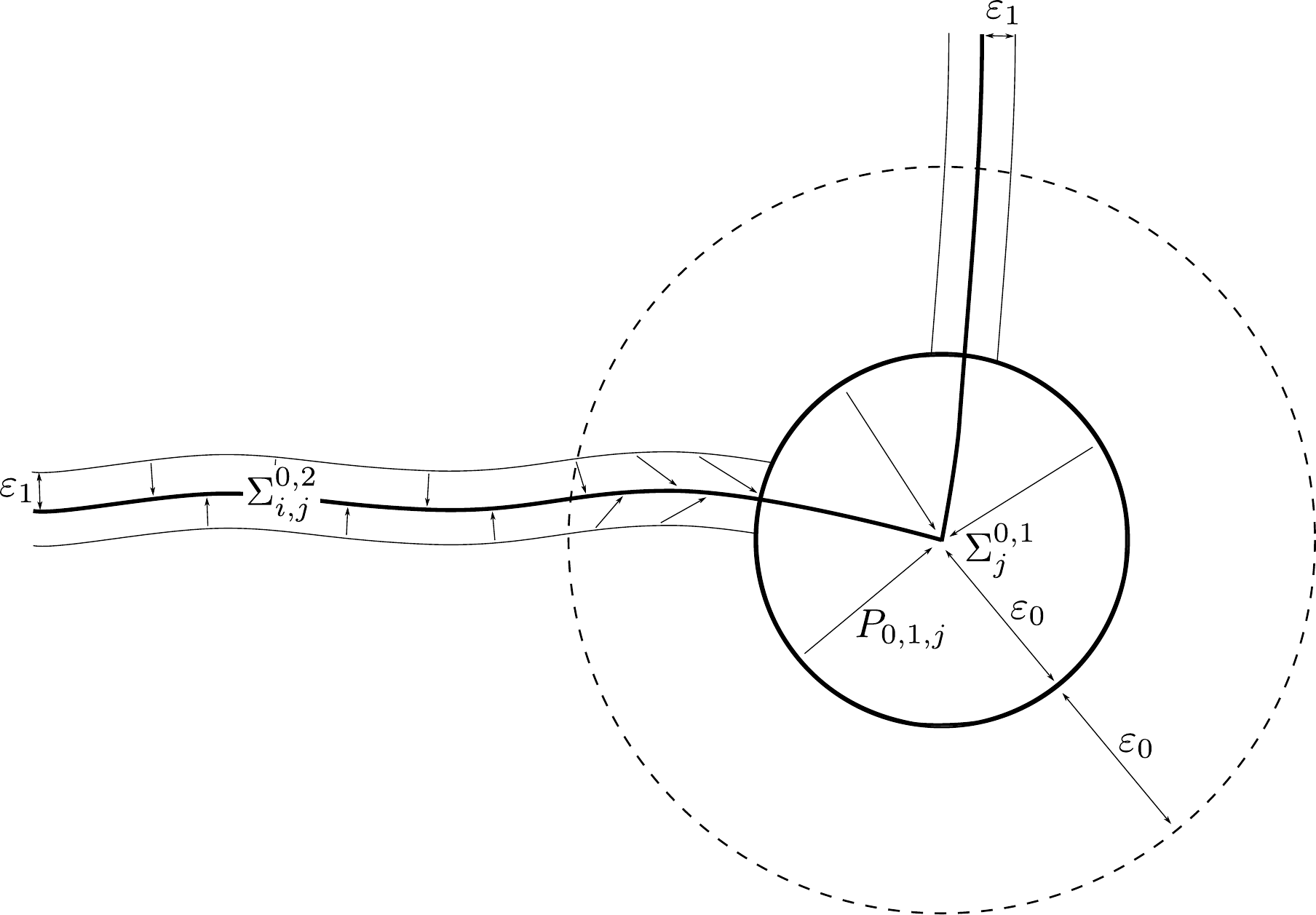}
\caption{A figure illustrating the construction of the transition maps at the intersection of different strata. In this case $m=3$ and $4\pi\left[\left(1+\alpha_{i}\right)+\left(1+\alpha_{j}\right)\right]<\rho<4\pi$ for any choice of  the indices $i, j$ such that $1\leq i < j \leq 3$. The space $\Sigma_{\rho,\underline{\alpha}}$ is made of three arcs joining the vertices $\delta_{p_{1}}, \delta_{p_{2}}, \delta_{p_{3}}$ in $C^{1}\left(\Sigma,g\right)^{\ast}$. Here we zoom around a vertex, say $\Sigma^{0,1}_{j}$ for some $j\in\left\{1,2,3\right\}$ and two arcs emanating from $\delta_{p_{j}}$ that correspond to two strata of dimension 1.}
\label{fig:circles}
\end{figure}

\begin{pfn}
Let $n$ denote the maximal dimension of an admissible stratum in $\Sigma_{\rho,\underline{\alpha}}$ and observe that obviously for any $j\leq n$ there exists only a finite number of strata having dimension $j$. After this preliminary remark, we define some numbers
\begin{displaymath}
\epsilon_{d_{l}}\ll\epsilon_{d_{l-1}}\ll\ldots\ll\epsilon_{d_{1}}\ll\epsilon_{d_{0}}\ll 1
\end{displaymath}
(with $n=d_{l}>d_{l-1}>\ldots>d_{1}>d_{0}=0$ denoting the dimensions of admissible strata of $\Sigma_{\rho,\underline{\alpha}}$) as follows. We choose $\epsilon$ so that for any admissible stratum of dimension 0 (say generically $\Sigma^{0,1}_{j}$) there is continuous projection from the $L^{1}(\Sigma)$ (normalized) functions in an $\epsilon$-neighborhood of  that $\Sigma^{0,1}_{j}$ onto $\Sigma^{0,1}_{j}$. This is possible by Lemma \ref{esistpr}. Then we consider all the strata of dimension $d_{1}$:  notice that it is not true in general that $d_{1}=1$ (see below for explicit examples), i. e. there could be dimensional gaps and in that case we just neglect those dimensions. However, we apply Lemma \ref{esistpr} again separately to each of these strata with $\epsilon_{0}=\frac{\epsilon}{4}$ and hence get a corresponding small $\widehat{\epsilon}$ and set $\epsilon_{d_{1}}=\frac{\widehat{\epsilon}}{4}$. We iterate the process and choose the numbers $\epsilon_{d_{2}},\ldots,\epsilon_{d_{l}}$ in the same way. \newline
For any $i\in \cup_{j=0}^{l}\left\{d_{j}\right\}$, let $f_{i}$ be a smooth non-increasing cut-off function such that
\begin{displaymath}
\left\{\begin{array}{ll}
f_{i}(t)=1,  & \textrm{for} \ t\leq \epsilon_{i};\\
f_{i}(t)=0,  & \textrm{for} \ t\geq 2\epsilon_{i}.
\end{array}\right.
\end{displaymath}
The next step consists in choosing the large number $\widehat{L}$, and this is essentially an elementary argument based on our concentration results above, Lemma \ref{l2sing}. The key point is that considering concentration at an appropriate scale, there exists a level $\widehat{L}$ such that for any $u\in H^{1}(\Sigma,g)$ with $J_{\rho,\underline{\alpha}}\left(u\right)\leq-\widehat{L}$ one has $d\left(\widetilde{h}e^{2u},\Sigma_{\rho,\underline{\alpha}}\right)<\epsilon_{d_{l}}$. Notice that here we are always assuming to work with functions normalized according to $\int_{\Sigma}\widetilde{h}e^{2u}\,dV_{g}=1$, which is no loss of generality since the functional is invariant under addition of constants to its argument. \newline
As a result, taken any $u\in H^{1}(\Sigma,g)$ with   $J_{\rho,\underline{\alpha}}\left(u\right)\leq-\widehat{L}$ there exists a smallest integer $j$ such that
$d\left(\widetilde{h}e^{2u},\Sigma^{k,l}_{\iota}\right)\leq \epsilon_{j}$ for some stratum $\Sigma^{k,l}_{\iota}$ in $\Sigma_{\rho,\underline{\alpha}}$ having dimension $j$. Hence, thanks to Lemma \ref{esistpr} and our choice of the $\epsilon_{i}$'s, the projection $P_{k,l,\iota}\left(\widetilde{h}e^{2u}\right)$ is well-defined and since (by definition of the index $j$) $d(\widetilde{h}e^{2u},\Sigma^{k',l'}_{\iota'})>\epsilon_{d}$ (where $d=\textrm{dim}\Sigma^{k',l'}_{\iota'}$) for any stratum $\Sigma^{k',l'}_{\iota'}$ such that $\Sigma^{k',l'}_{\iota'}\prec \Sigma^{k,l}_{\iota}$ the choice of such a stratum is unambiguous. Then we set
\begin{displaymath}
\Psi(u)=\odot_{\prec} U^{f_{d\left(k',l',\iota'\right)}\left(d\left(\widetilde{h}e^{2u},\Sigma^{k',l'}_{\iota'}\right)\right)}_{k',l',\iota'} \circ P_{k,l,\iota}\left(\widetilde{h}e^{2u}\right),
\end{displaymath}
where the symbol $\odot$ indicates a composition product which is extended to all homotopy operators $U^{t}_{\star}$ that correspond to strata $\Sigma^{k',l'}_{\iota'}\prec\Sigma^{k,l}_{\iota}$ and $d(k',l',\iota')$ is the dimension of the stratum $\Sigma^{k',l'}_{\iota'}$.

Notice that we are adopting the convention that the operators $U^{0}_{\star}$ that would in principle defined only locally are trivially extended to the whole $\Sigma_{\rho,\underline{\alpha}}$ as identity operator (this creates no problem because of property (i) in Lemma \ref{omotcoll}). The choice of extending the composition product to the strata $\Sigma^{k',l'}_{\iota'}\prec\Sigma^{k,l}_{\iota}$ is justified by Lemma \ref{inters}. The definition we have given depends in principle on the index $j$ which is a function of $u$. Nevertheless, since all distance functions from the strata are continuous and since $U^{1}_{\star}=P_{\star}$, this map $\Psi$ is actually well-defined and continuous in $u$.
\end{pfn}

The following property is a natural consequence of our construction.

\begin{cor}
\label{last}
Let $\Psi$ the projection map defined in the previous Lemma \ref{proi} and let $\widehat{L}\gg 1$ be the corresponding threshold value. If $\left(u_{n}\right)_{n\in\mathbb{N}}\subseteq J_{\rho,\underline{\alpha}}^{-\widehat{L}}$ and $\widetilde{h}e^{2u_{n}}\rightharpoonup \sigma$ for some $\sigma\in \Sigma_{\rho,\underline{\alpha}}$, then $\Psi(u_{n})\rightharpoonup \sigma$ in the weak sense.
\end{cor}

\

\

\section{Mapping $\Sigma_{\rho,\underline{\alpha}}$ into sublevels of $J_{\rho,\underline{\alpha}}$}

\medskip

In this section, we start by defining a very general class of bubbling functions  parameterized by the set $\Sigma_{\rho,\underline{\alpha}}$. Moreover, in order to perform a suitable min-max scheme in the proof of Theorem \ref{grande} (see Section 5), we want $J_{\rho,\underline{\alpha}}$ to attain arbitrarily negative values on such functions, this being true uniformly in $\sigma\in\Sigma_{\rho,\underline{\alpha}}$ when the scale parameter $\lambda$ tends to infinity. The difficult point in this step with respect to the regular case is that we need to take the presence of the singular points into account. By this reason, we introduce some sort of interpolation between the \textsl{regular} bubbling functions defined by \eqref{bubble} (more generally by \eqref{one}) and the singular bubbling functions defined by
\begin{equation}
\label{singb}
\varphi_{\alpha,\lambda,p}(y)=\log\left(\frac{\lambda^{1+\alpha}}{1+\left(\lambda d_{g}(p,y)\right)^{2(1+\alpha)}}\right),
\end{equation}
with $p=p_{j}$ for some $1\leq j\leq m$ and $\alpha=\alpha_{j}$ correspondingly.
For a small number $\delta>0$ we define the function $\gamma(\lambda,d)$ as
\begin{equation}
\label{gamma}
\gamma(\lambda,d)=\left\{\begin{array}{ll}
\alpha & \textrm{for} \ d<\delta\lambda^{-\frac{1}{1+\alpha}};\\
\gamma\in\left(0,\alpha\right) \ \textrm{s.t.} \ \lambda^{\frac{\gamma}{\alpha\left(1+\gamma\right)}}=\delta & \textrm{if} \ \delta\lambda^{-\frac{1}{1+\alpha}}<d\leq\delta;\\
 0 & \textrm{otherwise}.
\end{array}\right.
\end{equation}

Hence, for any $\sigma\in\Sigma_{\rho,\underline{\alpha}}$, say $\sigma=\sum_{i=1}^{k}t_{i}\delta_{x_{i}}$, we set
\begin{equation}
\label{multic}
\varphi_{\lambda,\sigma}(x)=\frac{1}{2}\log\left(\sum_{i}\frac{t_{i}\lambda^{2}}{\left(1+\lambda^{2}d_{g}(x,x_{i})^{2(1+\gamma_{i})}\right)^{2}}\right),
\end{equation}
where for any $i=1,\ldots,k$ we fix $\gamma_{i}=\gamma(\lambda,\min_{j}d_{g}(x_{i},p_{j}))$, and where the value $\alpha$ in \eqref{gamma} is the blow-up coefficient associated to the point $p_{j}$ realizing $\min_{j}d_{g}(x_{i},p_{j})$. To give sense to the definition \eqref{gamma} we must set $\alpha=0$ in case such a minimum is not smaller than $\delta$.

We are going to prove the following result.

\begin{pro}
\label{sprconti}
Let $\varphi_{\lambda,\sigma}$ be defined by \eqref{multic}. Then one has that
\begin{equation}
\label{lowen}
J_{\rho,\underline{\alpha}}\left(\varphi_{\lambda,\sigma}\right)\to-\infty \quad \textrm{as} \ \lambda\to+\infty,
\end{equation}
uniformly for $\sigma\in\Sigma_{\rho,\underline{\alpha}}$.
Moreover, there exists a universal constant $C>0$ (independent of $\lambda$) and coefficients $\widetilde{t}_{i}$ such that for any $i=1,\ldots,k$ 
\begin{displaymath}
\frac{t_{i}}{C}\leq \widetilde{t}_{i}\leq C t_{i}
\end{displaymath}
and
\begin{equation}
\label{ptifissi}
\widetilde{h}e^{2\varphi_{\lambda,\sigma}}\rightharpoonup \sum_{i=1}^{k}\widetilde{t}_{i}\delta_{x_{i}} \quad \textrm{as} \  \lambda\to+\infty.
\end{equation}
\end{pro}

In order to make the proof of this proposition more direct and effective, we choose to state the estimates for the Dirichlet energy term as a separate lemma, whose proof is postponed to the second part of this section.

\begin{lem}
\label{energy}
Let $\sigma=\sum_{i=1}^{n}t_{i}\delta_{x_{i}}$ and, correspondingly, $\mathcal{J}=\left\{x_{1},\ldots,x_{n}\right\}$. Then we have
\begin{equation}
\label{direner}
\int_{\Sigma}\left|\nabla\varphi_{\lambda,\sigma}\right|^{2}\,dV_{g}\leq 8\pi \chi\left(\mathcal{J}\right)\left(1+o_{\delta}(1)\right)\log\lambda+C_{\delta}.
\end{equation}
\end{lem}

We now prove Proposition \ref{sprconti}.\newline

\begin{pfn} Suppose some small number $\delta>0$ is fixed (the way to do this will be clear from the sequel).
We start by studying the integral $\int_{\Sigma}\varphi_{\lambda,\sigma}dV_{g}$. To this aim, notice that there exists a constant $C_{\delta}>0$ such that
\begin{displaymath}
-\log\lambda-C_{\delta}\leq \varphi_{\lambda,\sigma}(y)\leq \log\lambda \quad \textrm{in} \ \cup_{i=1}^{k}B_{\delta}(x_{i}),
\end{displaymath}
and
\begin{displaymath}
\left|\varphi_{\lambda,\sigma}(y)+\log\lambda\right|\leq C_{\delta} \quad \textrm{in} \ \Sigma\setminus\cup_{i=1}^{k}B_{\delta}(x_{i}).
\end{displaymath}
These estimates imply
\begin{equation}
\label{intbar}
\int_{\Sigma}\varphi_{\lambda,\sigma}\,dV_{g}=-(1+o_{\delta}(1))\log\lambda+O_{\delta}(1) \quad \textrm{as} \ \lambda\to+\infty.
\end{equation}
As our second step, we move to the study of the exponential term in the functional. We want to prove that
\begin{equation}
\label{estexp}
\log\int_{\Sigma}\widetilde{h}e^{2\varphi_{\lambda,\sigma}}\,dV_{g}=O(1)\quad \textrm{as} \ \lambda\to+\infty,
\end{equation}
more precisely we want to exhibit a constant $C$ such that
\begin{equation}
\label{Cindl}
\frac{1}{C}\leq\int_{\Sigma}\widetilde{h}(x)\frac{\lambda^{2}}{\left(1+\lambda^{2}d_{g}(x,x_{i})^{2(1+\gamma_{i})}\right)^{2}}\,dV_{g}\leq C,
\end{equation}
independently on $\lambda$ and for any possible value of the index $i$.
It should be clear that such a result also implies the second part of the thesis.
We need to split our manifold into three parts. First of all, it is clear that
\begin{equation}
\label{intexpext}
\int_{\Sigma\setminus B_{3\delta}(x_{i})}\widetilde{h}(x)\frac{\lambda^{2}}{\left(1+\lambda^{2}d_{g}(x,x_{i})^{2(1+\gamma_{i})}\right)^{2}}\,dV_{g}\leq \frac{C_{\delta}}{\lambda^{2}}.
\end{equation}
With respect to the other terms, it is necessary to consider two different cases, depending on whether $\min_{j=1,\ldots,m}d_{g}(x_{i},p_{j})\leq\delta$ or $\min_{j=1,\ldots,m}d_{g}(x_{i},p_{j})>\delta$. In the latter case, we can further divide the integral into $B_{\delta/2}(x_{i})$ and its complement with respect to $B_{3\delta}(x_{i}).$ In the second set the estimate is analogous to \eqref{intexpext}, while for the first set we do the computation in geodesic normal coordinates centered at $x_{i}\in\Sigma$.
In these coordinates one has
\begin{equation}
\label{propd}
dV_{g}=\left(1+o_{\delta}\left(1\right)\right)dx; \quad 1+\lambda^{2}d\left(x,x_{i}\right)^{2\left(1+\gamma_{i}\right)}=\left(1+o_{\delta}(1)\right)\left(1+\lambda^{2}\left|x-x_{i}\right|^{2(1+\gamma_{i})}\right),
\end{equation}
where we are implicitly identifying each point on the manifold $\Sigma$ (near $x_{i}$) with its normal coordinates. From \eqref{propd}, since in this case $\widetilde{h}$ is uniformly bounded from above and below by positive constants in $B_{\delta/2}(x_{i})$, one gets
\begin{equation}
\label{propd'}
\frac{1}{C_{\delta}}\leq\frac{\int_{B_{\delta/2}(x_{i})}\widetilde{h}(x)\frac{\lambda^{2}}{\left(1+\lambda^{2}d_{g}(x,x_{i})^{2(1+\gamma_{i})}\right)^{2}}\,dV_{g}}{\int_{\widehat{B}_{\delta/2}(x_{i})}\frac{\lambda^{2}}{\left(1+\lambda^{2}\left|x-x_{i}\right|^{2(1+\gamma_{i})}\right)^{2}}\,dx} \leq C_{\delta}.
\end{equation}
Here $\widehat{B}_{\delta/2}(x_{i})$ stands for a set in $\mathbb{R}^{2}$ that satisfies $B_{\left(1+o_{\delta}\left(1\right)\right)\frac{\delta}{2}}\left(x_{i}\right)\subseteq\widehat{B}_{\delta/2}(x_{i})\subseteq B_{\left(1+o_{\delta}\left(1\right)\right)\frac{\delta}{2}}\left(x_{i}\right)$. We are assuming $\min_{j=1,\ldots,m}d_{g}(x_{i},p_{j})>\delta$, so by \eqref{gamma} we simply have $\gamma_{i}=0$ and hence it is enough to consider the integral
\begin{displaymath}
\int_{\widehat{B}_{\delta/2}(x_{i})}\frac{\lambda^{2}}{\left(1+\lambda^{2}\left|x-x_{i}\right|^{2}\right)^{2}}\,dx.
\end{displaymath}
By a change of variables and elementary estimates, we conclude
\begin{displaymath}
\int_{\widehat{B}_{\delta/2}(x_{i})}\frac{\lambda^{2}}{\left(1+\lambda^{2}\left|x-x_{i}\right|^{2}\right)^{2}}\,dx=\int_{\widehat{B}_{\lambda\delta/2}(0)}\frac{1}{\left(1+\left|y\right|^{2}\right)^{2}}\,dy=C_{0}+O(\lambda^{-2}),
\end{displaymath}
being $C_{0}$ a fixed positive constant. As a result, in case $\min_{j=1,\ldots,m}d_{g}(x_{i},p_{j})>\delta$ we obtain \eqref{Cindl}.
Let us then turn to the harder case $\min_{j=1,\ldots,m}d_{g}(x_{i},p_{j})\leq\delta$. Here the singularities and their blow-up rate come into play.
Call $p$ the unique singular point that realizes $\min_{j=1,\ldots,m}d_{g}(x_{i},p_{j})$ and use geodesic coordinates centered at $p$. In these coordinates the approximation formulas \eqref{propd} still hold and so also \eqref{propd'} adapted to our case, hence
\begin{displaymath}
\frac{1}{C_{\delta}}\leq\frac{\int_{B_{3\delta}(x_{i})}\widetilde{h}(x)\frac{\lambda^{2}}{\left(1+\lambda^{2}d_{g}(x,x_{i})^{2(1+\gamma_{i})}\right)^{2}}\,dV_{g}}
{\int_{\widehat{B}_{3\delta}(x_{i})}\frac{\lambda^{2}\left|x\right|^{2\alpha}}{\left(1+\lambda^{2}\left|x-x_{i}\right|^{2(1+\gamma_{i})}\right)^{2}}\,dx}\leq C_{\delta}.
\end{displaymath}
Once again, we make the change of variables $y=\lambda^{\frac{1}{1+\gamma_{i}}}\left(x-x_{i}\right)$ and therefore
\begin{displaymath}
\int_{\widehat{B}_{3\delta}(x_{i})}\frac{\lambda^{2}\left|x\right|^{2\alpha}}{\left(1+\lambda^{2}\left|x-x_{i}\right|^{2(1+\gamma_{i})}\right)^{2}}\,dx
\end{displaymath}
\begin{displaymath}
=\int_{\widehat{B}_{3\delta\lambda^{\frac{1}{1+\gamma_{i}}}}(0)}\left|\lambda^{-\frac{1}{1+\gamma_{i}}}y+x_{i}\right|^{2\alpha}\frac{\lambda^{2}}{\left(1+\left|y\right|^{2(1+\gamma_{i})}\right)^{2}}\lambda^{-\frac{2}{1+\gamma_{i}}}\,dy
\end{displaymath}
\begin{equation}
\label{coeffcompl}
=\int_{\widehat{B}_{3\delta\lambda^{\frac{1}{1+\gamma_{i}}}}(0)}\frac{\left|\lambda^{\frac{\gamma_{i}-\alpha}{\left(1+\gamma_{i}\right)\alpha}}y+\lambda^{\frac{\gamma_{i}}{\left(1+\gamma_{i}\right)\alpha}}x_{i}\right|^{2\alpha}}{\left(1+\left|y\right|^{2(1+\gamma_{i})}\right)^{2}}\,dy.
\end{equation}
Now, we need to study this integral according to the different possible alternatives given by definition \eqref{gamma}. If we are in the first alternative of the definition of $\gamma_{i}$, the last integral becomes
\begin{equation}
\label{final}
\int_{\widehat{B}_{3\delta\lambda^{\frac{1}{1+\gamma_{i}}}}(0)}\frac{\left|y+v\right|^{2\alpha}}{\left(1+\left|y\right|^{2\left(1+\alpha\right)}\right)^{2}}\,dy
\end{equation}
where $v$ is a vector in $\mathbb{R}^{2}$ whose norm is uniformly bounded in $\lambda$ by some constant, say $C$. Since clearly $\delta\lambda^{\frac{1}{1+\gamma_{i}}}\to+\infty$ for $\lambda\to+\infty$, we can assume $\lambda$ so big that $\delta\lambda^{\frac{1}{1+\gamma_{i}}}\geq 2C$ and so the previous integral \eqref{final} is surely bounded from below. On the other hand, the same integral is less than the integral over $\mathbb{R}^{2}$ of the same function, which is uniformly bounded from above since the decay of the integrand at infinity is of order $\left|y\right|^{-4-2\alpha}$ and we are working with $\alpha\in\left(-1,0\right)$. So, if this alternative occurs we get \eqref{Cindl}.
\newline
In the second alternative for the definition of $\gamma_{i}$, the scalar $\lambda^{\frac{\gamma_{i}}{\left(1+\gamma_{i}\right)\alpha}}$ is exactly equal to $\delta$ and the coefficient of $y$ in \eqref{coeffcompl} is uniformly bounded. Hence, to get a lower bound, it is enough to integrate over a ball of radius $\delta^{2}$, while for an upper bound we mimic the previous argument, since the decay rate is $2\alpha-4-4\gamma_{i}<-2$ and the coefficient is uniformly bounded. This completes the proof of \eqref{Cindl}.
\newline
Now we just need to put together the previous estimates with the results claimed in Lemma \ref{energy} Indeed, combining \eqref{intbar}, \eqref{estexp} and \eqref{direner}, we find the uniform estimate
\begin{displaymath}
J_{\rho,\underline{\alpha}}(\varphi_{\lambda,\sigma})\leq\left(8\pi\chi\left(\mathcal{J}\right)-2\rho\right)\left(1+o_{\delta}(1)\right)\log\lambda+C_{\delta}
\end{displaymath}
and assuming $\delta$ is chosen sufficiently small this implies the thesis \eqref{lowen}.
\end{pfn} \newline

Let us go back to the proof of Lemma \ref{energy}. \newline

\begin{pfn}
To avoid too tedious notation we denote simply by $\varphi$ the function $\varphi_{\lambda,\sigma}$. We have:
\begin{displaymath}
\nabla\varphi(x)=\frac{1}{2}\frac{\sum_{i}\frac{-2t_{i}\lambda^{2}\left[1+\lambda^{2}d_{g}(x,x_{i})^{2\left(1+\gamma_{i}\right)}\right]\lambda^{2}\left(1+\gamma_{i}\right)d_{g}\left(x,x_{i}\right)^{2\gamma_{i}}\nabla_{g}d_{g}\left(x,x_{i}\right)^{2}}{\left[1+\lambda^{2}d_{g}\left(x,x_{i}\right)^{2\left(1+\gamma_{i}\right)}\right]^{4}}}{\sum_{i}\frac{t_{i}\lambda^{2}}{\left[1+\lambda^{2}d_{g}\left(x,x_{i}\right)^{2\left(1+\gamma_{i}\right)}\right]^{2}}},
\end{displaymath}
and so since the function $d_{g}\left(\bullet,x_{i}\right)$ is 1-Lipschitz this implies
\begin{displaymath}
\left|\nabla\varphi\left(x\right)\right|\leq\frac{\sum_{i}\frac{2t_{i}\left(1+\gamma_{i}\right)\lambda^{2}d_{g}\left(x,x_{i}\right)^{2\gamma_{i}+1}}{\left[1+\lambda^{2}d_{g}\left(x,x_{i}\right)^{2\left(1+\gamma_{i}\right)}\right]^{3}}}{\sum_{i}\frac{t_{i}}{\left[1+\lambda^{2}d_{g}\left(x,x_{i}\right)^{2\left(1+\gamma_{i}\right)}\right]^{2}}}.
\end{displaymath}
Via the following basic manipulation
\begin{displaymath}
\lambda^{2}d_{g}\left(x,x_{i}\right)^{2\gamma_{i}+1}=\lambda^{\frac{1}{1+\gamma_{i}}}\left[\lambda^{2}d_{g}(x,x_{i})^{2\left(1+\gamma_{i}\right)}\right]^{\frac{2\gamma_{i}+1}{2\left(\gamma_{i}+1\right)}}
\end{displaymath}
\begin{displaymath}
\leq\lambda^{\frac{1}{1+\gamma_{i}}}\left[1+\lambda^{2}d_{g}(x,x_{i})^{2\left(1+\gamma_{i}\right)}\right]^{\frac{2\gamma_{i}+1}{2\left(\gamma_{i}+1\right)}}
\end{displaymath}
we then obtain
\begin{equation}
\label{nfm}
\left|\nabla\varphi\left(x\right)\right|\leq\frac{\sum_{i}\frac{2t_{i}\left(1+\gamma_{i}\right)\lambda^{\frac{1}{1+\gamma_{i}}}}{\left[1+\lambda^{2}d_{g}\left(x,x_{i}\right)^{2\left(1+\gamma_{i}\right)}\right]^{2+\frac{1}{2\left(1+\gamma_{i}\right)}}}}{\sum_{i}\frac{t_{i}}{\left[1+\lambda^{2}d_{g}\left(x,x_{i}\right)^{2\left(1+\gamma_{i}\right)}\right]^{2}}}\leq m(x)
\end{equation}
provided we define
\begin{displaymath}
m(x)=\max_{i=1,\ldots,k}\left\{\frac{2\left(1+\gamma_{i}\right)\lambda^{\frac{1}{1+\gamma_{i}}}}{\left[1+\lambda^{2}d_{g}\left(x,x_{i}\right)^{2\left(1+\gamma_{i}\right)}\right]^{\frac{1}{2\left(1+\gamma_{i}\right)}}}\right\}.
\end{displaymath}
Let us restrict ourselves to the case when there is only one singularity $p$ with weight $\alpha$, since this does not really affect the generality of the argument. \newline
After choosing a sufficiently large constant $C>0$ we can divide the manifold $\Sigma$ into the following sets:
\begin{displaymath}
\mathcal{A}=\cup_{i}B_{C\lambda^{-\frac{1}{1+\gamma_{i}}}}\left(x_{i}\right)=:\cup_{i}\mathcal{A}_{i}, \quad \quad \mathcal{B}=\Sigma\setminus{\mathcal{A}}.
\end{displaymath}
We start studying the function $m(x)$ on the set $\mathcal{B}$: first of all we have the inequality
\begin{equation}
\label{mxe}
m(x)\leq\max_{i}\left\{\frac{2\left(1+\gamma_{i}\right)}{d_{g}\left(x,x_{i}\right)}\right\}.
\end{equation}
Then, choose one point (say $x_{\overline{i}}$) for which the distance from $p$ is the smallest among the $x_{i}$'s. For any other index $j\neq \overline{i}$ and a (sufficiently small) $\delta>0$ we consider the sets
\begin{displaymath}
\mathcal{B}_{j}=\mathcal{B}\cap \left\{x: \frac{1+\gamma_{j}}{d_{g}(x,x_{j})}>\left(1+\delta\right)\frac{1+\gamma_{\overline{i}}}{d_{g}\left(x,x_{\overline{i}}\right)} \quad \textrm{and} \quad \frac{1+\gamma_{j}}{d_{g}(x,x_{j})}>\max_{k\neq\overline{i}}\frac{1+\gamma_{k}}{d_{g}\left(x,x_{k}\right)}\right\}.
\end{displaymath}
In $\mathcal{B}\setminus \cup_{j\neq\overline{i}}\mathcal{B}_{j}$ we have
\begin{displaymath}
\max_{i}\left\{\frac{1+\gamma_{i}}{d_{g}(x,x_{i})}\right\}\leq\left(1+\delta\right)\frac{1+\gamma_{\overline{i}}}{d_{g}\left(x,x_{\overline{i}}\right)}
\end{displaymath}
and so we can substitute this into \eqref{mxe} to get
\begin{equation}
\label{intbmbj}
\int_{\mathcal{B}\setminus\cup_{j\neq\overline{i}}\mathcal{B}_{j}}\left(m(x)\right)^{2}\,dV_{g}\leq
4\left(1+\delta\right)^{2}\int_{\mathcal{B}\setminus\cup_{j\neq\overline{i}}\mathcal{B}_{j}}\frac{\left(1+\gamma_{\overline{i}}\right)^{2}}{d_{g}\left(x,x_{\overline{i}}\right)^{2}}\,dV_{g}
\end{equation}
\begin{displaymath}
\leq 4\left(1+\delta\right)^{2}\int_{\mathcal{B}\setminus{B_{C\lambda^{-\frac{1}{1+\gamma_{\overline{i}}}}}\left(x_{\overline{i}}\right)}}\frac{\left(1+\gamma_{\overline{i}}\right)^{2}}{d_{g}\left(x,x_{\overline{i}}\right)^{2}}\,dV_{g}
\leq 8\pi\left(1+\delta\right)^{2}\left(1+\alpha\right)\log\lambda+C_{\delta}.
\end{displaymath}
In $\mathcal{B}_{j}$ we first need to observe that the following two inequalities hold:
\begin{displaymath}
\frac{1+\gamma_{j}}{d_{g}\left(x,x_{j}\right)}>\left(1+\delta\right)\frac{1+\gamma_{\overline{i}}}{d_{g}\left(x,x_{\overline{i}}\right)}\geq\left(1+\delta\right)\frac{1+\gamma_{j}}{d_{g}\left(x,x_{\overline{i}}\right)},
\end{displaymath}
since $\gamma_{\overline{i}}$ is the biggest among the $\gamma$'s because $x_{\overline{i}}$ is the closest point to the singularity $p$. This implies
\begin{equation}
\label{cfrconi}
d_{g}\left(x,x_{\overline{i}}\right)>\left(1+\delta\right)d_{g}\left(x,x_{j}\right) \ \textrm{in} \ \mathcal{B}_{j}.
\end{equation}
We need to examine in more detail what are the points that satisfy this inequality and this is done geometrically comparing graphs of different distance functions in $\Sigma\times\mathbb{R}$ that are respectively centered at $x_{\overline{i}}$ with slope $1$ and centered at $x_{j}$ with slope $\left(1+\delta\right).$ It is clear that there exists a constant $C_{\delta}$ such that the points verifying \eqref{cfrconi} are contained in the ball $B_{C_{\delta}d_{g}\left(x_{j},x_{\overline{i}}\right)}\left(x_{j}\right)$. Hence, just exploiting the definition of $\mathcal{B}_{j}$ we find that
\begin{equation}
\label{intbj}
\int_{\mathcal{B}_{j}}\left(m(x)\right)^{2}\,dV_{g}\leq 4\left(1+\delta\right)^{2}\int_{\mathcal{B}_{j}}\frac{\left(1+\gamma_{j}\right)^{2}}{d_{g}\left(x,x_{j}\right)^{2}}\,dV_{g}
\end{equation}
\begin{displaymath}
\leq 4\left(1+\delta\right)^{2}\int_{B_{C_{\delta}d_{g}\left(x_{j},x_{\overline{i}}\right)}\left(x_{j}\right)\setminus B_{C\lambda^{-\frac{1}{1+\gamma_{j}}}}\left(x_{j}\right)}\frac{\left(1+\gamma_{j}\right)^{2}}{d_{g}\left(x,x_{j}\right)^{2}}\,dV_{g}
\end{displaymath}
\begin{displaymath}
\leq 8\pi\left(1+\gamma_{j}\right)^{2}\left[\frac{1}{1+\gamma_{j}}\log\lambda-\log\frac{1}{d_{g}\left(x_{j},x_{\overline{i}}\right)}\right]\left(1+o_{\delta}\left(1\right)\right)+C_{\delta}.
\end{displaymath}

From the triangle inequality, we have that
\begin{displaymath}
d_{g}(x_{j},x_{\overline{i}})\leq d_{g}(x_{j},p)+d_{g}(p,x_{\overline{i}})\leq 2d_{g}(p,x_{j})
\end{displaymath}
and so via substitution in \eqref{intbj}
\begin{displaymath}
\int_{\mathcal{B}_{j}}\left(m\left(x\right)\right)^{2}\,dV_{g}\leq 8\pi(1+\gamma_{j})^{2}\left[\frac{1}{1+\gamma_{j}}\log\lambda-\log\frac{1}{d_{g}(p,x_{j})}\right]\left(1+o_{\delta}(1)\right)+C_{\delta};
\end{displaymath}
therefore, recalling the definition \eqref{gamma} $\frac{1}{d_{g}(p,x_{j})}\geq C^{-1}\lambda^{\frac{\gamma_{j}}{\left(1+\gamma_{j}\right)\alpha}}$ we conclude that
\begin{displaymath}
\frac{1}{\left(1+\delta\right)^{2}}\int_{\mathcal{B}_{j}}\left(m(x)\right)^{2}\,dV_{g}
\end{displaymath}
\begin{displaymath}
\leq 8\pi\left(1+\gamma_{j}\right)^{2}\left[\frac{1}{1+\gamma_{j}}\log\lambda-\frac{\gamma_{j}}{\left(1+\gamma_{j}\right)\alpha}\log\lambda\right]\left(1+o_{\delta}(1)\right)+C_{\delta}
\end{displaymath}
\begin{displaymath}
=8\pi\left(1+\gamma_{j}\right)\left(1-\frac{\gamma_{j}}{\alpha}\right)\log\lambda\left(1+o_{\delta}(1)\right)+C_{\delta}
\leq 8\pi\log\lambda\left(1+o_{\delta}(1)\right)+C_{\delta}.
\end{displaymath}
Lastly, putting together \eqref{nfm}, \eqref{intbmbj} and \eqref{intbj} we obtain
\begin{equation}
\label{intb}
\int_{\mathcal{B}}\left|\nabla_{\lambda,\sigma}\right|^{2}\,dV_{g}\leq 8\pi\left(k+\alpha\right)\left(1+o_{\delta}(1)\right)\log\lambda+C_{\delta}.
\end{equation}
As a second step, we have to study $\int_{\mathcal{A}}\left(m\left(x\right)\right)^{2}\,dV_{g}$.
We introduce new functions $f_{i}(x)$ that come into play because of the following inequality
\begin{displaymath}
\frac{\lambda^{\frac{1}{1+\gamma_{i}}}}{\left[1+\lambda^{2}d_{g}\left(x,x_{i}\right)^{2\left(1+\gamma_{i}\right)}\right]^{\frac{1}{2\left(1+\gamma_{i}\right)}}}\leq
C\frac{\lambda^{\frac{1}{1+\gamma_{i}}}}{1+\lambda^{\frac{1}{1+\gamma_{i}}}d_{g}\left(x,x_{i}\right)}
\end{displaymath}
\begin{displaymath}
\leq C\frac{1}{\lambda^{-\frac{1}{1+\gamma_{i}}}+d_{g}\left(x,x_{i}\right)}=:f_{i}(x).
\end{displaymath}
Fixing $x\in\mathcal{A}$ we want to maximize (or better find upper bounds for) $f_{i}(x)$ with respect to the index $i$. \newline
We consider first the case of $x$ belonging to $\mathcal{A}\cap B_{2\delta\lambda^{-\frac{1}{1+\alpha}}}\left(p\right)$. For $x_{i}$ also in $B_{2\delta\lambda^{-\frac{1}{1+\alpha}}}\left(p\right)$ the function $f_{i}$ is bounded by $C\lambda^{\frac{1}{1+\alpha}}$. Let us assume that $x_{i}$ lies outside $B_{2\delta\lambda^{-\frac{1}{1+\alpha}}}\left(p\right)$ instead: in this case
\begin{displaymath}
1+\gamma_{i}=\frac{\log\lambda}{\log\lambda+\alpha\left(\log\left|x_{i}\right|-\log\delta\right)} \ \Rightarrow \ \frac{1}{1+\gamma_{i}}=1+\alpha\frac{\log\left|x_{i}\right|-\log\delta}{\log\lambda}.
\end{displaymath}
This implies
\begin{displaymath}
\label{expr}
\frac{1}{\lambda^{-\frac{1}{1+\gamma_{i}}}+\left|x-x_{i}\right|}=\frac{1}{C_{\delta}\lambda^{-1}\left|x_{i}\right|^{-\alpha}+\left|x_{i}-x\right|}.
\end{displaymath}
Notice that in the last two equations we are working in geodesic normal coordinates and again identifying points on $\Sigma$ and their coordinates on the tangent space $T_{p}\Sigma$. To get an upper bound for the latter quantity, we have to estimate the infimum of $C_{\delta}\lambda^{-1}\left|x_{i}\right|^{-\alpha}+\left|x_{i}-x\right|$ for $\left|x_{i}\right|\geq 2\delta\lambda^{-\frac{1}{1+\alpha}}$. By trivial geometric arguments one finds that this is of order $\lambda^{-\frac{1}{1+\alpha}}$ and therefore by all these estimates we get that
\begin{displaymath}
\sup_{B_{2\delta\lambda^{-\frac{1}{1+\alpha}}}(p)}m(x)\leq \sup_{i} \sup_{B_{2\delta\lambda^{-\frac{1}{1+\alpha}}}\left(p\right)} f_{i}(x)\leq C\lambda^{\frac{1}{1+\alpha}}.
\end{displaymath}
As a result
\begin{equation}
\label{intmx2AB}
\int_{\mathcal{A}\cap B_{2\delta\lambda^{-\frac{1}{1+\alpha}}}\left(p\right)}\left(m\left(x\right)\right)^{2}\,dV_{g}\leq C.
\end{equation}
We have next to consider the case in which $x\in\mathcal{A}\setminus B_{2\delta\lambda^{-\frac{1}{1+\alpha}}}$. For $x_{i}$ inside $B_{\delta\lambda^{-\frac{1}{1+\alpha}}}$, by \eqref{gamma} it is $\gamma_{i}=\alpha$ so that the denominator in $f_{i}(x)$ is bounded below by $\lambda^{-\frac{1}{1+\alpha}}$ and hence $f_{i}(x)$ is bounded by $\lambda^{\frac{1}{1+\alpha}}$. So we have reduced the problem to the case $x_{i}$ lies outside of $B_{\delta\lambda^{-\frac{1}{1+\alpha}}}$. We use again the expression \eqref{expr} that has to be maximized in terms of the position of $x_{i}$. The problem can be reduced to the one-dimensional case in which $x_{i}$ moves along the half-line emanating from $p$ towards $x$. By means of elementary calculus we find that
\begin{equation}
\label{supmx}
m\left(x\right)\leq C\lambda^{\frac{1}{1+\gamma\left(\lambda,\left|x\right|\right)}}, \quad \quad \textrm{for} \ x\in\mathcal{A}\setminus B_{2\delta\lambda^{-\frac{1}{1+\alpha}}}\left(p\right)
\end{equation}
and hence
\begin{displaymath}
\int_{\mathcal{A}_{j}\setminus B_{2\delta\lambda^{-\frac{1}{1+\alpha}}}\left(p\right)}\left(m\left(x\right)\right)^{2}\,dV_{g}\leq C\lambda^{-\frac{2}{1+\gamma_{j}}}\sup_{x\in\mathcal{A}_{j}\setminus B_{2\delta\lambda^{-\frac{1}{1+\alpha}}}\left(p\right)}\lambda^{\frac{2}{1+\gamma\left(\lambda,\left|x\right|\right)}}.
\end{displaymath}

Recalling the definition of $\mathcal{A}_{j}$ and
\begin{displaymath}
\frac{1}{1+\gamma\left(\lambda,\left|x\right|\right)}=1+\min\left\{0,\alpha\frac{\log\left|x\right|-\log\delta}{\log\lambda}\right\},
\end{displaymath}
we find that for $\left|x_{j}\right|\geq\delta\lambda^{-\frac{1}{1+\alpha}}$
\begin{displaymath}
\int_{\mathcal{A}_{j}\setminus B_{2\delta\lambda^{-\frac{1}{1+\alpha}}}\left(p\right)}m\left(x\right)^{2}\,dV_{g}\leq C\left(1+C\frac{1}{\lambda^{\frac{1}{1+\gamma_{j}}}\left|x_{j}\right|}\right)^{\alpha}
\end{displaymath}
\begin{displaymath}
\leq C\left(1+C_{\alpha,\delta}\lambda^{-1}\left|x_{j}\right|^{-1-\alpha}\right)\leq C_{\alpha,\delta},
\end{displaymath}
while for $\left|x_{j}\right|\leq\delta\lambda^{-\frac{1}{1+\alpha}}$
\begin{displaymath}
\int_{\mathcal{A}_{j}\setminus B_{2\delta\lambda^{-\frac{1}{1+\alpha}}}\left(p\right)}m\left(x\right)^{2}\,dV_{g}\leq
C\lambda^{-\frac{1}{1+\alpha}}\lambda^{\frac{1}{1+\gamma\left(\lambda,2C\lambda^{-\frac{1}{1+\alpha}}\right)}}\leq C_{\alpha,\delta}.
\end{displaymath}
From the last two inequalities and \eqref{intmx2AB} we finally obtain
\begin{equation}
\label{inta}
\int_{\mathcal{A}_{j}\setminus B_{2\delta\lambda^{-\frac{1}{1+\alpha}}}\left(p\right)}\left|\nabla\varphi_{\lambda,\sigma}\right|^{2}\,dV_{g}\leq C_{\alpha,\delta}.
\end{equation}
Combining \eqref{inta} and \eqref{intb} we get
\begin{displaymath}
\int_{\Sigma}\left|\nabla\varphi_{\lambda,\sigma}\right|^{2}\,dV_{g}\leq 8\pi\left(k+\alpha\right)\left(1+o_{\delta}(1)\right)\log\lambda+C_{\alpha,\delta}.
\end{displaymath}
In the general case, i.e. when we deal with any number of singularities, the same argument works just with minor modifications and leads to \eqref{direner}.
\end{pfn}


Now, we have all the tools needed to go back to the previous section and show that the map $\Psi$ is topologically non-trivial, so that it is not homotopically equivalent to a constant. Actually, we show something more.

\begin{lem}
\label{nonprov}
If $\Phi_{\lambda}\left(\sigma\right)=\varphi_{\lambda,\sigma}$ according to formula \eqref{multic}, then for $\lambda$ sufficiently large the map $\sigma\to\left(\Psi\circ \Phi_{\lambda}\right)\left(\sigma\right)=\Psi(\varphi_{\lambda,\sigma})$ is homotopic to the identity in $\Sigma_{\rho,\underline{\alpha}}$. As a result, if (and only if) such space is not contractible the projection $\Psi$ is non-trivial.
\end{lem}

\begin{pfn} We know by Lemma \ref{sprconti} (see especially formula \eqref{ptifissi}) and the previous Corollary \ref{last}, that for any $\sigma\in\Sigma_{\rho,\underline{\alpha}}$, say $\sigma=\sum_{i} c_{i}\delta_{z_{i}}$, $\Psi(\varphi_{\lambda,\sigma})\rightharpoonup \widetilde{\sigma}=\sum_{i}\widetilde{c}_{i}\delta_{z_{i}}$ for $\lambda\to+\infty$. It is clear that the coefficients $\widetilde{c}_{i}$ depend continuously on $\sigma$ and so we can define the map
\begin{displaymath}
\Omega:\Sigma_{\rho,\underline{\alpha}}\to\Sigma_{\rho,\underline{\alpha}}, \  \Omega(\sigma)=\widehat{\sigma}.
\end{displaymath}
We observe that $\Omega$ is homotopically equivalent to the identity $Id$ in $\Sigma_{\rho,\underline{\alpha}}$ by means of the homotopy
\begin{displaymath}
\left(\sigma,t\right)\longmapsto \left(1-t\right)\Omega\left(\sigma\right)+t\sigma.
\end{displaymath}
Notice that this is well-defined because $\sigma$ and $\Omega(\sigma)$ only differ by the coefficients, but not on the centers of the Dirac masses (this was proved in Lemma \ref{sprconti}).
Moreover, by the very definition of $\Omega$, we know that  for $\lambda$ sufficiently large the composition map $\Psi\circ\Phi_{\lambda}$ is homotopic to $\Omega$ itself in $\Sigma_{\rho,\underline{\alpha}}$. By composition of these two homotopic equivalences we finally get that for large $\lambda$'s  $\Psi\circ\Phi_{\lambda}$ is homotopic to the identity on $\Sigma_{\rho,\underline{\alpha}}$, which is exactly what we had to prove.
\end{pfn}


\

\section{Existence of solutions}

\medskip

The tools presented in the previous sections are all we need to prove our main result, namely Theorem \ref{grande}, which is essentially an existence theorem for non-critical values of $\rho$ (depending on $\underline{\alpha}$), related with the number in the denominator of \eqref{eq:cl}.

\noindent Our plan is to use a general min-max scheme in the form of a suitable \textsl{topological cone} construction.\newline

\textbf{1.} \ \textbf{Min-max scheme}. We assume a threshold value $L\gg1$ is chosen according to Lemma \ref{proi} and, correspondingly, $\lambda$ is fixed so that the operator $\Phi_{\lambda}$ takes values in the sublevel $J_{\rho,\underline{\alpha}}^{-2L}$, this being possible thanks to Lemma \ref{sprconti}. In order to simplify our notation we will omit explicit dependence on $\lambda$ in the sequel. We define the topological cone over $\Sigma_{\rho,\underline{\alpha}}$ as follows:
\begin{displaymath}
\Theta_{\rho,\underline{\alpha}}=\left(\Sigma_{\rho,\underline{\alpha}}\times \left[0, 1\right]\right)/\left(\Sigma_{\rho,\underline{\alpha}}\times\left\{1\right\}\right),
\end{displaymath} where we are identifying all the points in $\Sigma_{\rho,\underline{\alpha}}\times\left\{1\right\}.$ Consequently, we consider the family of continuous maps
\begin{displaymath}
\mathcal{H}_{\rho,\underline{\alpha}}=\left\{\mathfrak{h}:\Theta_{\rho,\underline{\alpha}}\to H^{1}(\Sigma,g)\ : \ \mathfrak{h}(\sigma)=\varphi_{\sigma} \  \textrm{for every} \ \sigma\in\Sigma_{\rho,\underline{\alpha}}\right\},
\end{displaymath}
and then the number
\begin{displaymath}
\overline{\mathcal{H}}_{\rho,\underline{\alpha}}=\inf_{\mathfrak{h}\in\mathcal{H}}\sup_{\sigma\in\Theta}J_{\rho,\underline{\alpha}}(\mathfrak{h}(\sigma)).
\end{displaymath} We claim that under the assumption of Theorem \ref{grande} one has $\overline{\mathcal{H}}_{\rho,\underline{\alpha}}\geq-L$. It is worth proving first that the class $\mathcal{H}_{\rho,\underline{\alpha}}$ is not empty. To this aim, notice that the map
$\underline{\mathfrak{h}}\left(\sigma,t\right)=\left(1-t\right)\varphi_{\sigma}, \ \left(\sigma,t\right)\in \Sigma_{\rho,\underline{\alpha}}$
belongs to $\mathcal{H}_{\rho,\underline{\alpha}}$.\newline
Concerning the lower bound on the min-max value, we just need to argue by contradiction.
If it were $\overline{\mathcal{H}}_{\rho,\underline{\alpha}}<-L,$ then there should be a map $\mathfrak{h}$ such that its image $\mathfrak{h}(\Theta_{\rho,\underline{\alpha}})$ (which is a topological cone in $H^{1}(\Sigma,g)$) would be in $J_{\rho,\underline{\alpha}}^{-L}.$ As a consequence, the composite map
\begin{displaymath}
t\to \Psi\left(\mathfrak{h}\left(\sigma,t\right)\right), \ \sigma\in\Sigma_{\rho,\underline{\alpha}}
\end{displaymath}
would be a homotopy equivalence between $\Psi\left(\mathfrak{h}\left(0,\sigma\right)\right)=\Psi\circ\Phi\left(\sigma\right)$ and a constant map. On the other hand, we know that the function $\Psi\circ\Phi\left(\sigma\right)$ is homotopic to the identity in $\Sigma_{\rho,\underline{\alpha}}$ (see Lemma \ref{nonprov}) and hence, by composition the space $\Sigma_{\rho,\underline{\alpha}}$ would be contractible, a contradiction. Hence we deduce $\overline{\mathcal{H}}_{\rho,\underline{\alpha}}\geq-L$.
\newline

\medskip

\textbf{2.} \ \textbf{Existence on a dense set.} The scheme outlined in the previous step immediately leads to existence for a dense set of $\rho$'s (in a suitable neighborhood of a fixed value). This relies on a monotonicity trick by Struwe \cite{s} and exploited also in \cite{djlw2}.
\newline

\medskip

\textbf{3.} \ \textbf{Conclusion via blow-up analysis.}
Let us now deal with any $\rho\in\mathbb{R}\setminus\mathfrak{S}$ to conclude our existence argument. The basic idea is very simple: build a sequence of approximating values $\left(\rho_{n}\right)_{n\in\mathbb{N}}$ such that $\rho_{n}\to\rho$ and $\rho_{n}\in\Lambda, \ \forall \ n\in\mathbb{N}$. This is clearly possible because $\Lambda$ has full measure.
Due to Step 2 we find a sequence $(v_{n})_{n\in\mathbb{N}}$ of solutions of $(\ref{mod})_{\rho_{n}}$ and recalling the substitution performed in the introduction, we can build a corresponding sequence $(u_{n})_{n\in\mathbb{N}}$, where $u_{n}=v_{n}+\sum_{j=1}^{m}\alpha_{j}G_{p_{j}}$, such that for any $n\in\mathbb{N}$ the function $u_{n}$ solves Problem \eqref{sing} for the parameter $\rho_{n}$ (the parameters $\underline{\alpha}$ are assumed to be fixed). Hence, we just need some compactness result and possibly also some regularity argument.
But before coming to the main results of this section, let us spend few words on the regularity of such solutions $v_{n}$  and $u_{n}$. Let $(v,u)$ denote any of the couples $(v_{n},u_{n}), n\in\mathbb{N}$ where $v=u-\sum_{j=1}^{m}\alpha_{j}G_{p_{j}}$. Thanks to the Moser-Trudinger inequality and to the fact that by assumption $\alpha_{j}>-1$ for all $j'$s, one easily finds that there is an $r>1$ such that $\widetilde{h}e^{2v}\in L^{r}(\Sigma,g)$ and so, by help of standard elliptic estimates we get $v\in W^{2,r}(\Sigma,g)$ and hence by the Sobolev embedding this gives $v\in C^{\alpha}(\Sigma,g)$ for some $\alpha\in\left(0,1\right)$. Moreover, by applying these arguments on domains of $\Sigma$ bounded away from $\left\{p_{1},\ldots,p_{m}\right\}$ we find that $v\in C^{\infty}(\Sigma\setminus\left\{p_{1},\ldots,p_{m}\right\})$. However, it should be clear that we cannot hope such maximal regularity on all of our manifold $\Sigma$. As a result, $u$ is a smooth function far from the singularities and has blow-up points at the singularities that are completely described by the corresponding Green functions, so $u\simeq \log d_{g}(x,p_{j})^{\alpha_{j}}$ near $p_{j}$ since $v$ is a H\"older function on the whole $\Sigma$. Hence we might say that $v$ is the \textsl{regular part}, while $\sum_{j=1}^{m}\alpha_{j}G_{p_{j}}$ is the \textsl{singular part} of $u$, a solution of \eqref{sing}.
We now come to the study of the limit phenomena that occur for the sequence $v_{n}$ when $n\to\infty$.

\medskip

\begin{thm}[\cite{bt1}]
\label{altern}
Let $w_{n}$ be any sequence of solutions of problem $(\ref{mod})_{\rho_{n}}$ in $H^{1}(\Sigma,g)$ for values $\rho_{n}$ of the parameter with $\rho_{n}\to\overline{\rho}$ and such that there exists a constant $C$ with
\begin{displaymath}
\int_{\Sigma}\widetilde{h}e^{2w_{n}}\,dV_{g}\leq C, \quad \forall n\in \mathbb{N}.
\end{displaymath}
There exists a subsequence $(w_{n_{k}})_{k\in\mathbb{N}}$ for which the following alternative holds:\newline
\underline{either} $w_{n_{k}}$ is uniformly bounded in $L^{\infty}(\Sigma,g)$;\newline
\underline{or} $\max_{\Sigma}\left(w_{n_{k}}-\frac{1}{2}\log \int_{\Sigma}\widetilde{h}e^{2w_{n_{k}}}\right)\to+\infty$,\newline
and there exists a finite (blow-up) set $S=\left\{z_{1},\ldots,z_{l}\right\}\subseteq \Sigma$ such that:
\begin{enumerate}
\item{for any $j\in\left\{1,\ldots,l\right\}$, there exists a sequence $\left(x_{j,k}\right)_{k\in\mathbb{N}}$ such that $x_{j,k}\to z_{j},$ $w_{n_{k}}(x_{j,k})\to+\infty$ and $w_{n_{k}}\to-\infty$ uniformly on any compact set $K\subseteq\Sigma\setminus S$,}
\item{\begin{displaymath}
\rho_{n_{k}}\frac{\widetilde{h}e^{2w_{n_{k}}}}{\int_{\Sigma}\widetilde{h}e^{2w_{n_{k}}}\,dV_{g}}\rightharpoonup \sum_{j=1}^{l}\beta_{j}\delta_{z_{j}} \quad \textrm{in \ the \ sense \ of \ measures},
\end{displaymath}
with $\beta_{j}=4\pi$ for $z_{j}\notin \left\{p_{1},\ldots,p_{m}\right\},$ or $\beta_{j}=4\pi(1+\alpha_{i})$ in case $x_{j}=p_{i}$ for some $i\in\left\{1,\ldots,m\right\}$.}
\end{enumerate}

As a result, if this second alternative occurs, then $\rho\in\mathfrak{S}$ (defined by means of formula \eqref{defsv}).
\end{thm}

\begin{rem}
It should be noticed that this kind of result was first obtained by the authors of \cite{bt1} under the assumption $\alpha_{j}>0$ for every $j=1,\ldots,m$, but their argument works also in case the same parameters are negative. However, this requires some modifications, that are described in \cite{bmt}.
\end{rem}

This immediately gives what we need.

\begin{cor}
Assume $w_{n}$ is any family of solutions of $\eqref{mod}_{\rho}$ corresponding to values of $\rho$ belonging to a compact subset of $\mathbb{R}\setminus\mathfrak{S}$. Then $(w_{n})_{n\in\mathbb{N}}$ is uniformly bounded from above on $\Sigma$.
\end{cor}

More generally, we have the following

\begin{cor}[Concentration/Compactness]
Let $w_{n}$ be a sequence of solutions of $\eqref{mod}_{\rho}$. Then $w_{n}$ admits a subsequence that satisfies the following alternative:\newline
\underline{either} $w_{n}$ is uniformly bounded from above on $\Sigma$ and converges uniformly in $C^{\gamma}(\Sigma,g)$ for any $\gamma\in\left[0,\gamma_{0}\right)$ with $\gamma_{0}\in\left(0,1\right)$,\newline
\underline{or} the second case in Theorem \ref{altern} holds.
\end{cor}

This corollary is proved with no effort starting from Theorem \ref{altern}: in fact, if the first case occurs there, the extracted subsequence $(w_{n_{k}})_{k\in\mathbb{N}}$ is bounded and so the term $e^{2w_{n_{k}}}$ is also uniformly bounded in $L^{\infty}(\Sigma,g)$. The desired conclusion comes from a bootstrap argument and standard elliptic estimates.

\

\section{Examples and open problems}

\medskip

\

As outlined in the introduction of this article, Theorem \ref{grande} reduces the analytical problem of existence for equation \eqref{sing} to a purely topological problem. Basically, we are led to study the spaces $\Sigma_{\rho,\underline{\alpha}}$ for all admissible values of the parameters $\rho$ and $\underline{\alpha}$ or at least to determine whether or not they are contractible. When $m$ and $\rho$ are sufficiently large answering this question is definitely not trivial and indeed this is still an open problem.
In this section we first want to describe some applications of Theorem \ref{grande} and, as a result, we need to exhibit some specific cases of non-contractibility of the space $\Sigma_{\rho,\underline{\alpha}}$. This is primarily intended in order to give a visual and intuitive idea of the topological structure of such a space in some simple examples. We determine the labels of the singular points $p_{1},\ldots,p_{m}$ so that $\alpha_{1}\leq\alpha_{2}\leq\ldots\leq\alpha_{m}$ and, moreover, we always implicitly assume $\rho<8\pi$ and $\rho<4\pi\left(2+\alpha_{1}\right)$. Notice that we will repeatedly make use of the simple but enlightening Lemma \ref{inters} concerning the intersections of different strata of $\Sigma_{\rho,\underline{\alpha}}.$ \newline

\textbf{$k$-points configurations}.
Assume that $m\geq 1$ and the parameters $\rho$, $\underline{\alpha}$ satisfy the algebraic system
\begin{displaymath}
\left\{\begin{array}{ll}
\rho> 4\pi\left(1+\alpha_{i}\right),  & \textrm{for} \ 1\leq i\leq k;\\
\rho< 4\pi\left(1+\alpha_{i}\right),  & \textrm{for} \ k+1\leq i \leq m;\\
\rho< 4\pi\left[\left(1+\alpha_{i}\right)+\left(1+\alpha_{j}\right)\right],  & \textrm{for any couple of indices such that} \ 1\leq i<j\leq m,
\end{array}\right.
\end{displaymath}
for some integer $k$ such that $1\leq k\leq m$ with the convention that if $k=m$ this means $\alpha_{m+1}=0$ i.e.
\begin{displaymath}
\left\{\begin{array}{ll}
\rho> 4\pi\left(1+\alpha_{i}\right),  & \textrm{for} \ i\leq m;\\
\rho< 4\pi, \\
\rho< 4\pi\left[\left(1+\alpha_{i}\right)+\left(1+\alpha_{j}\right)\right],  & \textrm{for any couple of indices such that} \ 1\leq i<j\leq m.
\end{array}\right.
\end{displaymath}
In these cases the space $\Sigma_{\rho,\underline{\alpha}}$ simply consists of $k$ points, indeed
\begin{displaymath}
\Sigma_{\rho,\underline{\alpha}}=\cup_{i\leq k}\left\{\delta_{p_{i}}\right\}.
\end{displaymath}
This means that the very low sublevels of the functional $J_{\rho,\underline{\alpha}}$ mirror this topology in the sense that they have $k$ (arc-wise) connected components, each one being contractible. \newline As a consequence, if $\Sigma_{\rho,\underline{\alpha}}$ only consists of strata having dimension 0, then this space is contractible if and only if $k=1$.\newline

\textbf{Graphs with loops} Following a naive ordering by increasing topological complexity, immediately after $k$-points configurations we find graphs. It is well known and easy to prove that a (finite) connected graph is contractible if (and only if) it does not contain loops. Observe that by Lemma \ref{inters} the nodes of our graphs are the (admissible ones among) vertices $\delta_{p_{j}}, j\in\left\{1,2,\ldots,m\right\}$ and the edges are the $1$-simplices corresponding to strata $\Sigma^{0,2}_{\star}$. The case $m\leq 2$ is trivial and so assume $m\geq 3$: if we exclude the presence of strata of dimension greater or equal than 2, to get a loop we just need to require that there exists a triplet of pairwise distinct indices, say $\left\{i_{1},i_{2},i_{3}\right\}\subseteq\left\{1,\ldots,m\right\}$ such that $\rho>4\pi\left[\left(1+\alpha_{i_{j}}\right)+\left(1+\alpha_{i_{l}}\right)\right]$ for any choice of $i_{j}\neq i_{l}$. But since we are always assuming the ordering $\alpha_{1}\leq\alpha_{2}\leq\ldots\leq\alpha_{m}$ we have proved the following:
\begin{thm}
Assume the space of formal barycenters $\Sigma_{\rho,\underline{\alpha}}$ only consists of strata having dimension 0 or 1. Then necessary and sufficient conditions for the non-contractibility of that space are given by:\newline

\underline{either}
\begin{displaymath}
\left\{\begin{array}{ll}
m\geq 2, \ 2\leq k \leq m;\\
\rho> 4\pi\left(1+\alpha_{i}\right),  & \textrm{for} \ 1\leq i\leq k;\\
\rho< 4\pi\left(1+\alpha_{i}\right),  & \textrm{for} \ k+1\leq i \leq m;\\
\rho< 4\pi\left[\left(1+\alpha_{i}\right)+\left(1+\alpha_{j}\right)\right],  & \textrm{for any couple of indices such that} \ 1\leq i<j\leq m,
\end{array}\right.
\end{displaymath}

\underline{or}
\begin{displaymath}
\left\{\begin{array}{ll}
m\geq 3;\\
\rho> 4\pi\left[\left(1+\alpha_{i}\right)+\left(1+\alpha_{j}\right)\right],  & \textrm{for any couple of indices such that} \ 1\leq i<j\leq 3.
\end{array}\right.
\end{displaymath}
\end{thm}

Observe that requiring that the space $\Sigma_{\rho,\underline{\alpha}}$ does not contain strata of dimension greater than two is obtained by means of the conditions $\rho<4\pi$ and $\rho<4\pi\sum_{i=1}^{3}\left(1+\alpha_{i}\right)$, the second one being necessary only if $m\geq 3.$\newline

\textbf{Linear handles} Let us go back to the case described in Section 1. Indeed, let us require $m=2$ and
\begin{displaymath}
\left\{\begin{array}{ll}
\rho> 4\pi,  \\
\rho> 4\pi\left[\left(1+\alpha_{1}\right)+\left(1+\alpha_{2}\right)\right].  
\end{array}\right.
\end{displaymath}
We may embed $\Sigma_{\rho,\underline{\alpha}}$ in $\mathbb{R}^{3}$ obtaining a compact surface with a one-dimensional handle, that is an arc joining the singular points $p_{1}$ and $p_{2}$. The topological non-triviality is clear and in fact can be proved by elementary tools. We can generalize this example by taking many linear handles instead of only one and this happens whenever $m\geq 3$ and the parameters satisfy the algebraic inequalities

\begin{displaymath}
\left\{\begin{array}{ll}
\rho> 4\pi,  \\
\rho>4\pi\left[\left(1+\alpha_{1}\right)+\left(1+\alpha_{2}\right)\right],\\
\rho< 4\pi\sum_{i=1}^{3}\left(1+\alpha_{i}\right).
\end{array}\right.
\end{displaymath}

\textbf{$2-$simplices over $\Sigma$}.

\begin{figure}[h]
\centering
\includegraphics[scale=0.6]{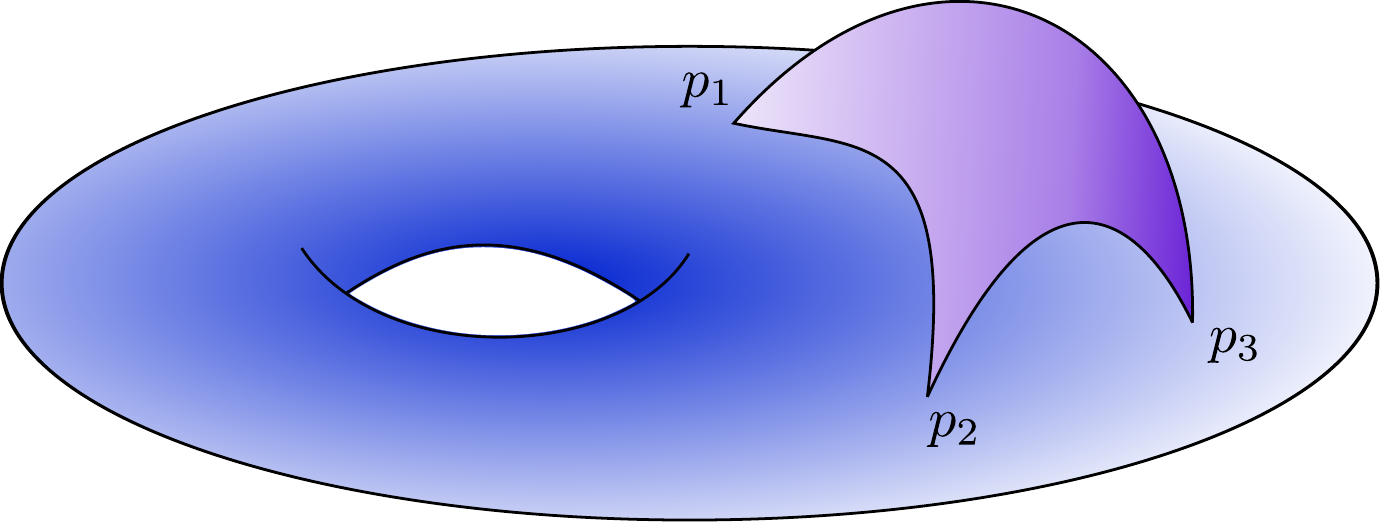}
\caption{A sketch of the space $\Sigma_{\rho,\underline{\alpha}}$ in case the parameters $\rho$ and $\underline{\alpha}$ satisfy the algebraic system \eqref{eqvela}. Notice the purple 2-dimensional \textsl{sail}.}
\label{fig:toro3}
\end{figure}

In case $m=3$ and
\begin{equation}
\label{eqvela}
\left\{\begin{array}{ll}
\rho> 4\pi,  \\
\rho>4\pi\sum_{i=1}^{3}\left(1+\alpha_{i}\right).
\end{array}\right.
\end{equation}
(recall we are always assuming $\rho\in\left(0,8\pi\right)$) we get that the space of formal barycenters $\Sigma_{\rho,\underline{\alpha}}$ is homeomorphic to the union (again via gluing at the singular points) of $\Sigma$ and a sort of \textsl{sail} (a 2-simplex). \newline

Indeed, the study of a wide range of special cases leads to formulate the following conjecture.

\begin{definition}
Given the parameters $\rho$ and $\underline{\alpha}$, we say that the corresponding model space $\Sigma_{\rho,\underline{\alpha}}$ is $p_{j}-$stable for some index $j\in{1,2,\ldots,m}$ if one of the following two equivalent conditions holds:
\begin{enumerate}
\item{Whenever $\sigma\in \Sigma_{\rho,\underline{\alpha}}$ then $(1-t)\sigma+t\delta_{p_{j}} \in \Sigma_{\rho,\underline{\alpha}} \ \forall\ t \in \left[0,1\right]$;}
\item{Whenever $k\in\mathbb{N}$ and a multi-index $\iota$ are such that \begin{displaymath}4\pi\left[k+\sum_{i\in\iota}\left(1+\alpha_{i}\right)\right]<\rho \end{displaymath} then also \begin{displaymath}4\pi\left[k+\sum_{i\in\left\{j\right\}\cup\iota}\left(1+\alpha_{i}\right)\right]<\rho.\end{displaymath}}
\end{enumerate}
\end{definition}

\begin{rem}
The condition given at point 2. of the previous definition cannot in general be simplified. Indeed, one could at first be lead to claim that $p_{j}-$stability is also equivalent to the much simpler requirement that if $\overline{k}, \overline{\iota}$ are such that 
\begin{equation}
\label{nns}
4\pi\left[\overline{k}+\sum_{i\in\overline{\iota}}\left(1+\alpha_{i}\right)\right]=\max_{adm}4\pi\left[k+\sum_{i\in\iota}\left(1+\alpha_{i}\right)\right]
\end{equation}
(where the maximum is taken over all admissible singular values, see \eqref{defsv}), then $j\in\iota$. In fact, this condition is necessary, but not sufficient for $p_{j}-$stability, as shown by the elementary example of the case $m=2,  4\pi(1+\alpha_{2})<\rho<4\pi\left[(1+\alpha_{1})+\left(1+\alpha_{2}\right)\right]$ for $j=2$.
\end{rem}

\begin{rem}
Notice that the corresponding notion of $q-$stability, for generic (namely regular) $q\in\Sigma$ would be meaningless since it is easily checked that $\Sigma_{\rho,\underline{\alpha}}$ is \textsl{never} $q-$stable for regular $q$. Notice also that indeed we can always reduce to consider the case $j=1$ by noticing that \textsl{if $\Sigma_{\rho,\underline{\alpha}}$ is $p_{j}-$stable for some index $j$, then it is necessarily $p_{1}-$stable}. To this aim, we argue as follows: suppose $k, \iota$ are given so that $4\pi\left[k+\sum_{i\in\iota}\left(1+\alpha_{i}\right)\right]<\rho$. There are two cases: either $j\in \iota$ or $j\notin\iota$. In the second alternative, the thesis is trivial since by assumption $\alpha_{1}\leq\alpha_{j}$. In the first, define the multi-index $\widetilde{\iota}$ by replacing in $\iota$ the index $j$ by the index $1$ (if $1\in\iota$, then we simply erase the index $j$). Clearly, $4\pi\left[k+\sum_{i\in\widetilde{\iota}}\left(1+\alpha_{i}\right)\right]<\rho$ and, thanks to the $p_{j}-$stability assumption we get $4\pi\left[k+\sum_{i\in{\left\{j\right\}}\cup\widetilde{\iota}}\left(1+\alpha_{i}\right)\right]<\rho$ which is equivalent to $4\pi\left[k+\sum_{i\in{\left\{1\right\}}\cup\iota}\left(1+\alpha_{i}\right)\right]<\rho$, so $\Sigma_{\rho,\underline{\alpha}}$ is $p_{1}-$stable.
\end{rem}

The reason why we are interested in $p_{1}-$stability is that if $\Sigma_{\rho,\underline{\alpha}}$ is $p_{1}-$stable, then it is contractible or, more precisely, it deformation-retracts onto $\delta_{p_{1}}$ in the ambient space $C^{1}\left(\Sigma,g\right)^{\ast}$ by means of the homotopy map $H:\Sigma_{\rho,\underline{\alpha}}\times\left[0,1\right]\rightarrow \Sigma_{\rho,\underline{\alpha}}$ given by $H(\sigma,t)=(1-t)\sigma+t\delta_{p_{1}}.$ It seems likely that the converse is also true:

\begin{conjecture}[topological version]
The space of formal barycenters $\Sigma_{\rho,\underline{\alpha}}$ is contractible if and only if it is $p_{1}-$stable.
\end{conjecture}

\begin{figure}[h]
\centering
\includegraphics[scale=0.7]{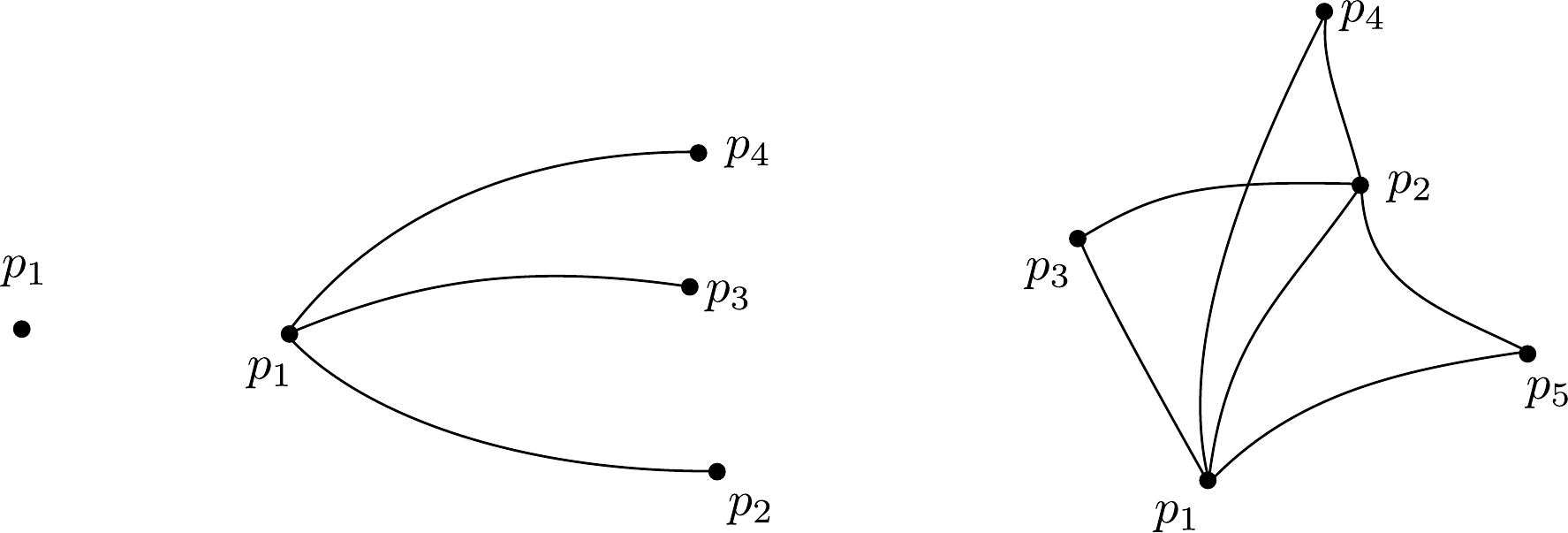}
\caption{Some prototypes of contractibility for the space $\Sigma_{\rho,\underline{\alpha}}$.}
\label{fig:grafo}
\end{figure}

\begin{example}
Let us describe the examples of Figure \ref{fig:grafo}.
\begin{enumerate}
\item[(i)]{In the first case $\Sigma_{\rho,\underline{\alpha}}$ is reduced to a single point,  $\Sigma_{\rho,\underline{\alpha}}=\Sigma^{0,1}_{1}$.}
\item[(ii)]{In the second case $\Sigma_{\rho,\underline{\alpha}}$ consists of three 1-simplices having a vertex in common, $\Sigma_{\rho,\underline{\alpha}}=\Sigma^{0,2}_{12}\cup \Sigma^{0,2}_{13}\cup \Sigma^{0,2}_{14}$.}
\item[(iii)]{In the third case $\Sigma_{\rho,\underline{\alpha}}$ consists of three 2-simplices having a 1-simplex in common, $\Sigma_{\rho,\underline{\alpha}}=\Sigma^{0,3}_{123}\cup\Sigma^{0,3}_{124}\cup\Sigma^{0,3}_{125}$.}
\end{enumerate}
\end{example}

Despite these examples, notice that we do not require all the strata belonging to a contractible $\Sigma_{\rho,\underline{\alpha}}$ to have the same dimension.

\vspace{2mm}

The previous conjecture can be immediately turned into algebraic form.

\begin{conjecture}[algebraic version]
\label{conalg}
The space of formal barycenters $\Sigma_{\rho,\underline{\alpha}}$ is  NOT contractible if and only if there exist a number $n\in\mathbb{N}$ and a set $\iota\subseteq\left\{2,3,\ldots,m\right\}$ such that $ card\left(\iota\right)\geq 1$ and 
\begin{displaymath}
\rho>4\pi\sum_{i\in \iota}\left(1+\alpha_{i}\right) \ \wedge \ \rho<4\pi\sum_{i\in\left\{1\right\}\cup \iota}\left(1+\alpha_{i}\right).
\end{displaymath}
\end{conjecture}

This is easily proved, by almost elementary methods, when we reduce to the case $\Sigma_{\rho,\underline{\alpha}}$ only consists of strata having dimension less than 3 or when $\rho<8\pi$ and in many other special cases, but a fair general proof seems to be rather hard. For instance, observe that when $4\pi<\rho<4\pi\left(2+\alpha_{1}\right)$ the thesis follows by simply considering the Mayer-Vietoris exact sequence in homology
\begin{equation*}
\begin{CD}
\ldots @>>> H_{2}(A\cap B;\mathbb{Z})  @>>>  H_{2}(A;\mathbb{Z})\oplus H_{2}(B;\mathbb{Z}) @>>> H_{2}(X;\mathbb{Z}) @>>> H_{1}(A\cap B;\mathbb{Z}) @>>> \ldots
\end{CD}
\end{equation*}
where $X=\Sigma_{\rho,\underline{\alpha}}$, $A$ is an $\varepsilon$-neighborhood of $\Sigma\hookrightarrow\Sigma_{\rho,\underline{\alpha}}$ and $B$ is an $\varepsilon$-neighborhood of $\Sigma_{\rho,\underline{\alpha}}\setminus\Sigma,$ for some small $\varepsilon$. Indeed, $A\cap B$ can be deformation-retracted onto a finite and non-empty set of points, hence $H_{1}(A\cap B;\mathbb{Z})=0$ and $H_{2}(A\cap B;\mathbb{Z})=0$ and therefore $H_{2}(X;\mathbb{Z})\approx H_{2}(A;\mathbb{Z})\oplus H_{2}(B;\mathbb{Z})$, this being non-trivial since $H_{2}(A;\mathbb{Z})\approx \mathbb{Z}$.
Anyway, in case this conjecture were true we could derive a large class of existence theorems directly by checking algebraic inequalities that involve the parameters $\rho$ and $\underline{\alpha}$.  \newline

Another related question naturally arises: Are the algebraic conditions above (Conjecture \ref{conalg}) only sufficient or also necessary for existence? It has recently been proved (see \cite{blt}) that in some cases of non-contractibility actually no solutions may exist. The class of tools that are used for this kind of argument are variations on the Pohozaev identity. So one could at first be led to claim that whenever $\Sigma_{\rho,\underline{\alpha}}$ is contractible we do not have existence. In fact, such converse implication seems rather unlikely. The reason is that even in very special cases (for instance $\Sigma=\mathbb{T}^{2}$ with the flat metric and $m=1$) it should be possible obtain solutions for Problem \eqref{sing} as local minima for the functional $J_{\rho,\underline{\alpha}}$ by means of a smart choice of the datum $h$. Similar techniques are often used in order to obtain multiplicity results, as shown for instance in \cite{st}, \cite{dmfr}, \cite{dmfr2} or \cite{dmfr3} and hence there is good reason to believe that in the next few years also this question will be answered in general situations. \newline

As a final remark, it should be highlighted that the definition of the space of formal barycenters $\Sigma_{\rho,\underline{\alpha}}$ given in Section 1 is believed to apply, without modifications, also to the more general case when the parameters $\alpha_{1},\ldots,\alpha_{m}$ are real numbers, some of which possibly being positive. This has already been proved in \cite{mr} in the case $\rho<8\pi$ and $0<\alpha_{i}\leq 1$. The general case is under our current investigation and, if verified, it would directly lead to a wide range of applications, primarily to the problem of prescribing Gaussian curvature for orbifolds with conical singularities, which we plan to specifically treat in a forthcoming paper.
\


\newpage

\small


\begin{thebibliography}{99}
\setcounter{footnote}{0}

\bibitem{aub}\textsc{T. Aubin}, \textit{Some Nonlinear Problems in Riemannian Geometry}, SMM Springer-Verlag, Berlin, 1998.
\bibitem{bc} \textsc{A. Bahri, J. M. Coron}, \textit{On a nonlinear elliptic equation
involving the critical Sobolev exponent: the effect of the
topology of the domain}, Comm. Pure Appl. Math. \textbf{41} (1988), 253-294.
\bibitem{bdm}\textsc{D. Bartolucci, F. De Marchis, A. Malchiodi}, \textit{Supercritical conformal metrics on surfaces with conical singularities}, Int. Math. Res. Not., to appear.
\bibitem{blt}\textsc{D. Bartolucci, C. S. Lin, G. Tarantello}, work in progress.
\bibitem{bmt}\textsc{D. Bartolucci, E. Montefusco}, \textit{Blow-up analysis, existence and qualitative properties of solutions of the two-dimensional Emden-Fowler equation with singular potential}, Math. Meth. Appl. Sci. \textbf{30} (2007), 2309-2327.
\bibitem{bt1}\textsc{D. Bartolucci, G. Tarantello}, \textit{Liouville type equations with singular data and their applications to periodic multivortices for the electroweak theory}, Comm. Math. Phys. \textbf{229} (2002), 3-47.
\bibitem{bt2}\textsc{D. Bartolucci, G. Tarantello}, \textit{The Liouville equation with singular data: a concentration-compactness principle via a local representation formula}, J. Diff. Eq. \textbf{185} (2002), 161-180.
\bibitem{bm}\textsc{H. Brezis, F. Merle}, \textit{Uniform estimates and blow up behaviour for solutions of $-\Delta u=V(x)e^{u}$ in two dimensions}, Comm. Part. Diff. Eq. \textbf{16} (1991), no. 8-9, 1223-1253.
\bibitem{clmp1}\textsc{E. Caglioti, P. L. Lions, C. Marchioro, M. Pulvirenti}, \textit{A special class of stationary flows for two dimensional Euler equations: A statistical mechanics description}, Comm. Math. Phys. \textbf{143} (1992), 501-525.
\bibitem{clmp2}\textsc{E. Caglioti, P. L. Lions, C. Marchioro, M. Pulvirenti}, \textit{A special class of stationary flows for two dimensional Euler equations: A statistical mechanics description, part II}, Comm. Math. Phys. \textbf{174} (1995), 229-260.
\bibitem{cm}\textsc{A. Carlotto, A. Malchiodi}, \textit{A class of existence results for the singular Liouville equation}, C.R.A.S. Serie Mathematique, \textbf{349} (2011), no. 3-4, 161-166.
\bibitem{cl1} \textsc{C. C. Chen, C. S. Lin}, \textit{Sharp estimates for solutions
of multi-bubbles in compact Riemann surfaces}, Comm. Pure Appl. Math. \textbf{55} (2002), 728-771.
\bibitem{cldeg} \textsc{C. C. Chen, C. S. Lin}, \textit{Topological degree for a mean field equation on Riemann surfaces}, Comm. Pure Appl. Math. \textbf{56} (2003), no. 12, 1667-1727.
\bibitem{clsing}   \textsc{C. C. Chen, C. S. Lin}, \textit{Mean field equations of Liouville type with singular data: sharper estimates}. Discrete Contin. Dyn. Syst. \textbf{28} (2010), no. 3, 1237-1272.
\bibitem{cl3} \textsc{C. C. Chen, C. S. Lin}, \textit{A degree counting formula for singular Liouville-type equation and its application to multi vortices in electroweak theory}, in preparation.
\bibitem{cli1}\textsc{W. X. Chen, C. Li}, \textit{Prescribing gaussian curvature on surfaces with conical singularities}, J. Geom. Anal. \textbf{1} (1991), no.4, 359-372.
\bibitem{cli2}\textsc{W. X. Chen, C. Li}, \textit{Classifications of solutions of some nonlinear elliptic equations}, Duke Math. J. \textbf{63} (1991), 615-622.
\bibitem{dmfr}\textsc{F. De Marchis}, \textit{Multiplicity result for a scalar field equation on compact surfaces}, Comm. in Part. Diff. Eq. \textbf{33} (2008), no. 12, 2208-2224.
\bibitem{dmfr2}\textsc{F. De Marchis}, \textit{Generic multiplicity for a scalar field equation on compact surfaces}, preprint, 2010.
\bibitem{dmfr3}\textsc{F. De Marchis}, \textit{Multiplicity of solutions for a mean field equation on compact surfaces}, preprint, 2010.
\bibitem{dem}\textsc{M. Del Pino, P. Esposito, M. Musso}, \textit{Two-dimensional Euler flows with concentrated vorticities}, preprint, 2008.
\bibitem{djlw2}\textsc{W. Ding, J. Jost, J. Li, G. Wang}, \textit{Existence results for mean field equations}, Ann. Inst. Henri Poincaré \textbf{16} (1999), 653-666.
\bibitem{zind}\textsc{Z. Djadli}, \textit{Existence result for the mean field problem on Riemann surfaces of all genuses}, Comm. Contemp. Math. \textbf{10} (2008), no. 2, 205-220.
\bibitem{dm}\textsc{Z. Djadli, A. Malchiodi}, \textit{Existence of conformal metrics with constant $Q-$curvature}, Ann. Math. \textbf{168} (2008), pp. 813-858.
\bibitem{esp}\textsc{P. Esposito}, \textit{Blow-up solutions for a Liouville equation with singular data}, SIAM J. Math. Anal. \textbf{36} (2005), no. 4, 1310-1345.
\bibitem{kw}\textsc{J. Kazdan, F. Warner}, \textit{Curvature functions for compact 2-manifolds}, Ann. Math. \textbf{99} (1974), 14-47.
\bibitem{kiess}\textsc{M. H. K. Kiessling}, \textit{Statistical mechanics of classical particles with logaritmic interaction}, Comm. Pure Appl. Math. \textbf{46} (1993), 27-56.
\bibitem{lai}\textsc{C. H. Lai (ed.)}, \textit{Selected Papers on Gauge Theory of Weak and Electromagnetic Interactions}, World Scientific, Singapore, 1981.
\bibitem{lp}\textsc{J. M. Lee, T. H. Parker}, \textit{The Yamabe problem}, Bull. of the A. M. S. \textbf{17} (1987), 37-91.
\bibitem{li}\textsc{Y. Y. Li}, \textit{Harnack type inequality: The method of moving planes}, Comm. Math. Phys. \textbf{200} (1999), no. 2, 421-444.
\bibitem{ls1}\textsc{Y. Y. Li, I. Shafrir}, \textit{Blow-up analysis for solutions of $-\Delta u=Ve^{u}$ in dimension two}, Indiana Univ. Math. J. \textbf{43} (1994), no. 4, 1255-1270.
\bibitem{luc}\textsc{M. Lucia}, \textit{A deformation lemma with an application to a mean field equation}, Topol. Methods Nonl. Anal. \textbf{30} (2007), no. 1, 113-138.
\bibitem{mal2}\textsc{A. Malchiodi}, \textit{Morse theory and a scalar field equation on compact surfaces}, Adv. Diff. Eq. \textbf{13} (2008), no. 11-12, 1109-1129.
\bibitem{mal1}\textsc{A. Malchiodi}, \textit{Topological methods for an elliptic equation with exponential nonlinearities}, Discr. Cont. Dyn. Syst. \textbf{21} (2008), no. 1, 277-294.
\bibitem{mr}\textsc{A. Malchiodi, D. Ruiz}, \textit{New improved Moser-Trudinger inequalities and singular Liouville equations on compact surfaces}, preprint, 2010.
\bibitem{mos}\textsc{J. Moser}, \textit{A sharp form of an inequality of N. Trudinger}, Indiana Univ. Math. J. \textbf{20} (1971), 1077-1092.
\bibitem{s}\textsc{M. Struwe}, \textit{The existence of surfaces of constant mean curvature with free boundaries}, Acta Math. \textbf{160} (1988), no. 1-2, 19-64.
\bibitem{st}\textsc{M. Struwe, G. Tarantello}, \textit{On multivortex solutions in Chern-Simons gauge theory}, Boll. Unione Mat. Ital., Sez. B Art. Ric. Mat. \textbf{8} (1998), 109-121.
\bibitem{tar1}\textsc{G. Tarantello}, \textit{Analytical aspects of Liouville-type equations with singular sources}, Handbook of Differential Equations. Stationary Partial Differential Equations, Elsevier-Sciences \textbf{1}, M.Chipot, P.Quittner Eds., 2006.
\bibitem{tar4}\textsc{G. Tarantello}, \textit{Self-Dual Gauge Field Vortices: An analytical approach}, PNLDE \textbf{72}, Birkh\"auser Boston Inc., Boston MA, 2007.
\bibitem{tr}\textsc{M. Troyanov},\textit{Prescribing curvature on compact surfaces with conical singularities}, Trans. A. M. S. \textbf{324} (1991), no. 2, 793-821.
\bibitem{luzh2} \textsc{L. Zhang}, \textit{Asymptotic behavior of blowup solutions for
elliptic equations with exponential nonlinearity and singular data}, Comm. Contemp.
Math. \textbf{11}, No. 3 (2009), 395-411.

\end{thebibliography}
\end{document}